\documentclass[12pt]{article}
\usepackage{amssymb}
\usepackage{amsmath}
\usepackage{amsthm}
\usepackage{comment}
\usepackage{graphicx}
\begin{document}

\def\bx{{\bf x}}
\def\ee{\varepsilon}
\def\myproof{\noindent {\bf Proof.}\\}
\def\pdf{PDF}

\newcommand{\removableFootnote}[1]{}

\newtheorem{theorem}{Theorem}
\newtheorem{conjecture}[theorem]{Conjecture}
\newtheorem{lemma}[theorem]{Lemma}



\title{Stochastically Perturbed Sliding Motion in Piecewise-Smooth Systems.}
\author{
D.J.W.~Simpson and R.~Kuske\thanks{
The authors are indebted to Mario di Bernardo for many useful discussions
regarding this work.
}\\
Department of Mathematics\\
University of British Columbia\\
Vancouver, BC\\
Canada}
\maketitle

\begin{abstract}

Sliding motion is evolution on a switching manifold of a 
discontinuous, piecewise-smooth system of ordinary differential equations.
In this paper we quantitatively study the effects of
small-amplitude, additive, white Gaussian noise on stable sliding motion.
For equations that are static in directions parallel to the switching manifold,
the distance of orbits from the switching manifold approaches a quasi-steady-state density.
From this density we calculate the mean and variance for the near sliding solution.
Numerical results of a relay control system
reveal that the noise may significantly affect the period
and amplitude of periodic solutions with sliding segments.

\end{abstract}

\section{Introduction}
\label{sec:INTRO}

Nonsmoothness and noise are two
features of a dynamical system that may be the cause of important qualitative behaviour.
Hybrid and piecewise-smooth systems are utilized in a wide variety of fields
to model phenomena that involve switching,
impacts or other nonsmooth elements \cite{DiBu08,VaSc00,LeNi04,BaVe01,ZhMo03,PuSu06}.
Recent studies have successively explained novel
behaviour that may occur in such systems.
For instance, so-called discontinuity maps have been developed
and analyzed in order to understand border-collision scenarios in
piecewise-smooth systems of ordinary differential equations \cite{DaNo00,FrNo00,DiBu01}.
Moreover, parameter uncertainty, background vibrations and other sources of noise are
ubiquitous in real systems.
Studies of stochastic differential equations have, for example,
led to an understanding of noise-induced dynamics
such as stochastic resonance and coherence resonance
in excitable systems \cite{BeGe06,LiGa04,PiKu97,GaHa98}.
However, investigations into systems that are both piecewise-smooth
and involve noise are relatively uncommon.

One-dimensional, piecewise-linear maps with noise have been the subject of
some isolated investigations \cite{Gr05,Wa98,ZhSh06,ZhHu10,Ti02}.
In \cite{GrHo05}, noise-induced transitions for a two-dimensional,
piecewise-smooth system of
ordinary differential equations are explained through an analysis of a one-dimensional,
piecewise-linear return map with additive noise.
Simple vibro-impacting systems have been analyzed with stochastic averaging \cite{DiIo04}.
Exact results are attainable for a classical, unforced, linear oscillator
with elastic impacts \cite{DiMe79,FoBr96}, whereas more complex scenarios
have been investigated asymptotically and numerically \cite{SrPa05,FeXu09b}.
Stochasticity in switched control systems is
particularly important in regards to
robustness of an output signal to noise \cite{ChEl05,FeZh06,MhEl05}.
Here Lyapunov functions are an invaluable mathematical tool
for determining stability because they do not necessitate
smoothness in the vector field.
Noise-induced oscillatory motion has been studied in
piecewise-linear systems for which local linearity makes
some key calculations tractable \cite{LiSc00,SiKu11}.

In this paper we consider piecewise-smooth, stochastic differential equations
for which the underlying deterministic dynamics are described by the ODE
\begin{equation}
\dot{\bx} = F_i(\bx) \;, \qquad \bx \in \Omega_i \;,
\label{eq:Filippov}
\end{equation}
where each $\Omega_i \subset \mathcal{D} \subset \mathbb{R}^N$
is open, nonempty and pairwise-disjoint, $\cup_i \overline{\Omega}_i = \mathcal{D}$,
and each $F_i : \overline{\Omega}_i \to \mathbb{R}^N$ is a smooth function.
Equation (\ref{eq:Filippov}) is a {\em Filippov system} \cite{Fi88}
and well-suited to model phenomena that alternate
between different dynamical regimes, such as
vibrating systems experiencing impacts or friction \cite{WiDe00,Br99,BlCz99,AwLa03,Ib09},
and switching in electrical circuits \cite{BaVe01,ZhMo03,Ts03}.

\begin{figure}[t!]
\begin{center}
\setlength{\unitlength}{1cm}
\begin{picture}(14.5,5.25)
\put(0,0){\includegraphics[height=5.25cm]{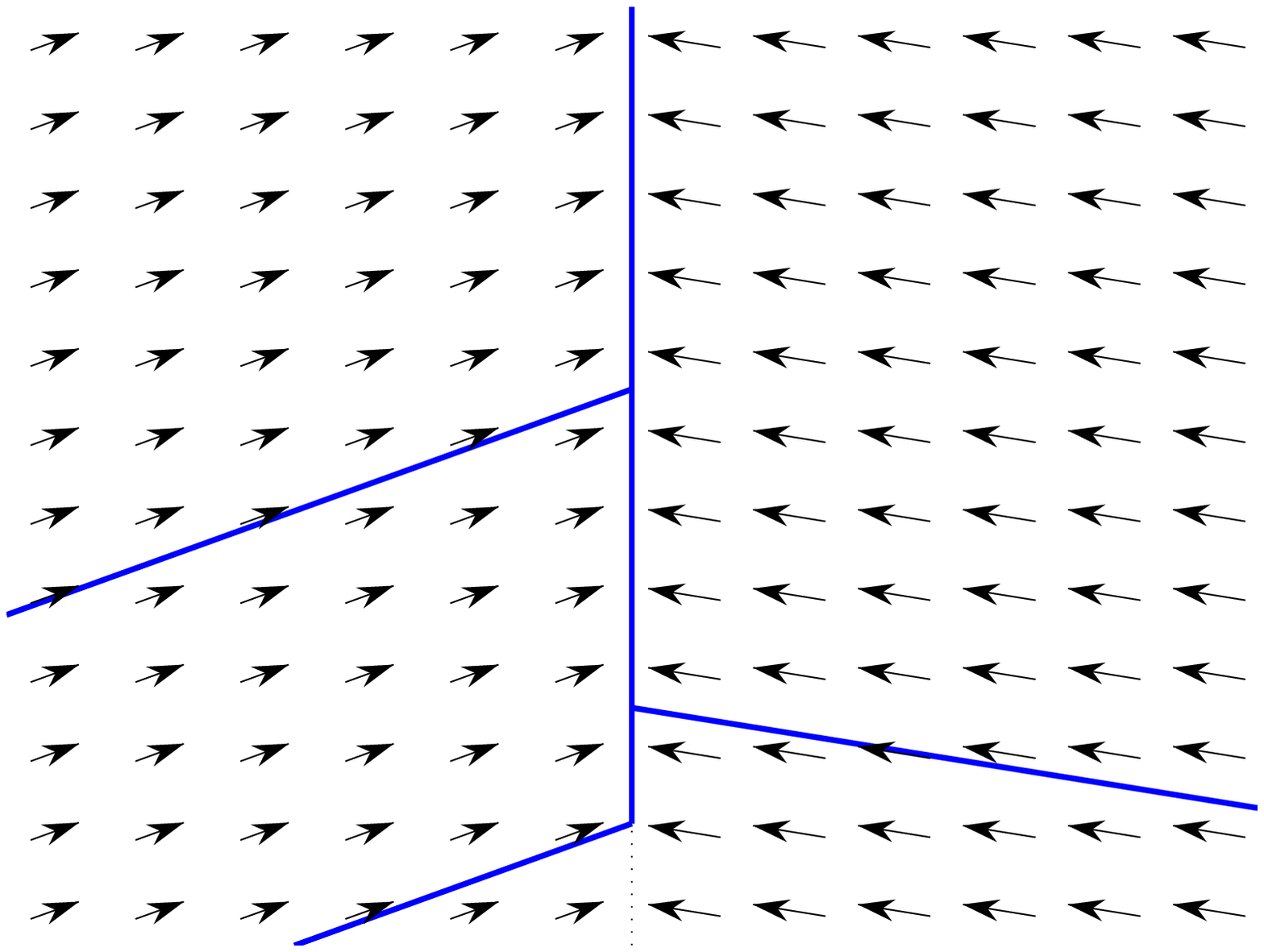}}
\put(7.5,0){\includegraphics[height=5.25cm]{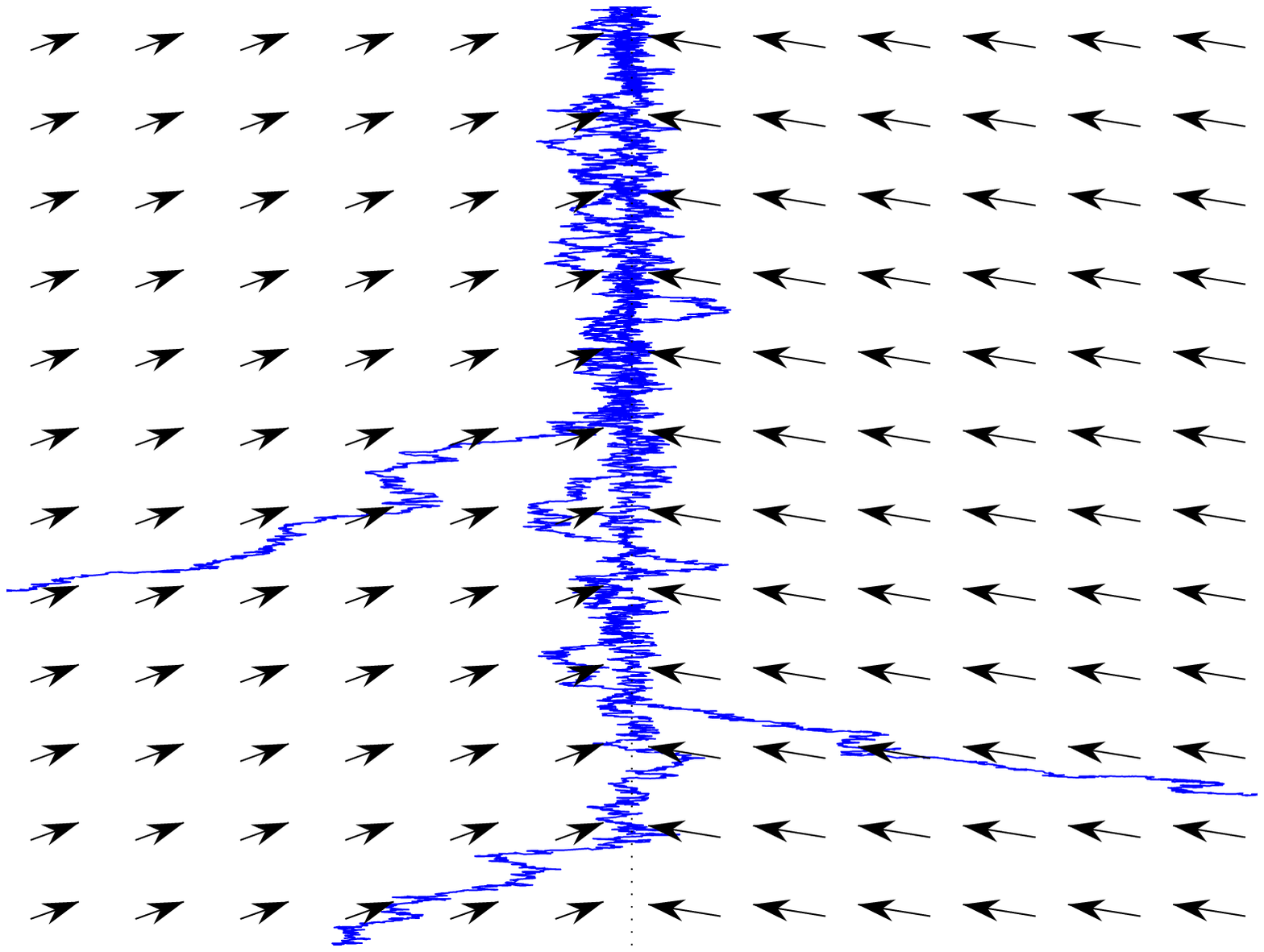}}
\put(1,5.25){\large \sf \bfseries A}
\put(8.5,5.25){\large \sf \bfseries B}
\end{picture}
\caption{Schematics of a Filippov system near a switching manifold that
attracts orbits from both sides in the absence of noise, panel A,
and with small amplitude additive noise, panel B.
\label{fig:noisyFlowExSlide}
}
\end{center}
\end{figure}

Boundaries between the neighbouring subdomains, $\Omega_i$, are codimension-one
surfaces termed {\em switching manifolds}.
Often, a section of a switching manifold
has the property that
on either side of the manifold the vector field points towards the manifold.
In this case any orbit that reaches the switching manifold
becomes trapped on the manifold for some time.
The resulting motion on the switching manifold is known as {\em sliding motion},
Fig.~\ref{fig:noisyFlowExSlide}-A.
Formally this is achieved by Filippov's method \cite{Fi60,Fi88,DiBu08}
which defines a vector field on the switching manifold
by the unique convex combination of the two limiting
vector fields on either side
that is tangent to the switching manifold.
Sliding motion corresponds to the sticking phase of
stick-slip oscillators \cite{CaGi06,LuGe06}
and the coalesced regime of a piecewise-linear
relay control system \cite{Jo03,DiJo01}.
The addition of small noise pushes orbits off the switching manifold,
but large excursions are curbed by the deterministic component
of the system, Fig.~\ref{fig:noisyFlowExSlide}-B.
Thus the motion is balanced by the competing
actions of noise and drift.
In this paper we explore these dynamics more carefully.

We here summarize the basic effect of adding small noise to a smooth system.
Below we compare this to
our results for the piecewise-smooth system (\ref{eq:Filippov}).
Consider an $N$-dimensional system
\begin{equation}
\dot{\bx} = \Phi(\bx) \;,
\label{eq:smoothODE}
\end{equation}
where $\Phi$ is a smooth function.
The addition of small amplitude, time-independent, white, Gaussian noise gives the
stochastic differential equation
\begin{equation}
d\bx = \Phi(\bx) \, dt + \sqrt{\ee} B(\bx) \, dW(t) \;,
\label{eq:smoothSDE}
\end{equation}
where $0 < \ee \ll 1$, $B(\bx)$ is an $N \times N$ matrix with a smooth
dependency on $\bx$, and $W(t)$ is a standard Brownian motion \cite{Sc80,Ga85}.
Let $p_\ee(\bx,t|\bx_0)$ denote the transitional probability density function (\pdf)
for the point $\bx(t)$, given $\bx(0) = \bx_0$.
A straight-forward expansion\removableFootnote{
Let us do this expansion for the one-dimensional stochastic differential equation:
\begin{equation}
dx = \phi(x) \, dt + \sqrt{\ee} B(x) \, dW(t) \;.
\end{equation}
In $N$-dimensions the process is messier
(i.e.~tensors) but otherwise essentially the same.
Write
\begin{equation}
x(t) = x_0(t) + \sqrt{\ee} x_1(t) + \ee x_2(t) + \ldots \;,
\end{equation}
where in practice I find it easier to let $\xi = \sqrt{\ee}$
and expand in $\xi$.
Then
\begin{eqnarray}
dx_0 + \sqrt{\ee} \, dx_1 + \ee \, dx_2 + \ldots &=&
\left( \phi(x_0) + \sqrt{\ee} \phi'(x_0) x_1 +
\ee \left( \frac{1}{2} \phi''(x_0) x_1^2 + \phi'(x_0) x_2 \right) + \ldots \right) \, dt
\nonumber \\
&&+~\left( \sqrt{\ee} B(x_0) + \ee B'(x_0) + \ldots \right) \, dW(t)
\end{eqnarray}
Collecting $O(1)$ terms shows that $x_0(t)$
equals the solution in the absence of noise.
Taking $O(\sqrt{\ee})$ terms gives
\begin{eqnarray}
dx_1 &=& \phi'(x_0) x_1 \, dt + b(x_0) \, dW(t) \;, \\
\Rightarrow~~x_1(t) &=&
\int_0^t {\rm e}^{\int_s^t \phi'(x_0(u)) \, du}
B(x_0(s)) \, dW(s) \;,
\end{eqnarray}
by using an integrating factor and substituting $x_1(0) = 0$.
This is a time-dependent Ornstein-Uhlenbeck process,
$x_1(t)$ is normally distributed and
\begin{equation}
\big< x_1(t) \big> = 0 \;, \qquad
{\rm Var}(x_1(t)) = \int_0^t {\rm e}^{2 \int_s^t \phi'(x_0(u)) \, du} B(x_0(s))^2 \, ds \;.
\end{equation}
Similarly taking $O(\ee)$ terms gives
\begin{eqnarray}
dx_2 &=& \left( \frac{1}{2} \phi''(x_0) x_1^2 + \phi'(x_0) x_2 \right) \, dt
+ B'(x_0) \, dW(t) \nonumber \\
\Rightarrow~~x_2(t) &=&
\frac{1}{2} \int_0^t \phi''(x_0(s)) x_1(s)^2
{\rm e}^{\int_s^t \phi'(x_0(u)) \, du} \, ds +
\int_0^t B'(x_0(s)) {\rm e}^{\int_s^t \phi'(x_0(u)) \, du} \, dW(s) \;,
\end{eqnarray}
using also $x_2(0) = 0$.
Consequently
\begin{equation}
\int_0^t \phi''(x_0(s)) \big< x_1(t)^2 \big>
{\rm e}^{\int_s^t \phi'(x_0(u)) \, du} \, ds \;.
\end{equation}
Therefore the difference between $\big< x(t) \big>$ and the
solution in the absence of noise is $\ee \big< x_2(t) \big> +
O \left( \ee^{\frac{3}{2}} \right)$.
For the lowest order term, $\big< x_2(t) \big>$,
to be nonzero we require nonlinearity in $\phi$ (i.e.~$\phi''(x) \ne 0$).
}
in powers of $\sqrt{\ee}$ reveals that the mean of $p_\ee(\bx,t|\bx_0)$
differs from the deterministic solution
(the solution to (\ref{eq:smoothODE}) with $\bx(0) = \bx_0$)
by $O(\ee)$, whereas deviations are $O(\sqrt{\ee})$ \cite{Ga85,GrVa99}.

In this paper we derive an analogous result for sliding motion
for which the above method of expansion does not work
because the vector field is discontinuous.
Instead we analyze a one-dimensional, discontinuous stochastic differential equation
for a quantity representing the distance from the switching manifold.
We find that the mean solution differs from
Filippov's deterministic sliding solution by $O(\ee)$,
and deviations are $O(\sqrt{\ee})$, matching the smooth case.
Moreover, the calculations suggest conditions necessary
for the noise to induce a change in the dynamics that is
not dominated by randomness.

Here we outline the remainder of the paper.
We motivate our work in \S\ref{sec:MOTIV} 
by illustrating that noise may significantly reduce the amplitude
and period of a solution to a prototypical relay control model
involving segments of sliding motion.
In \S\ref{sec:TOY} we introduce a
system of two-dimensional stochastic differential equations, (\ref{eq:e})-(\ref{eq:aLaRpos}),
that describes stochastically perturbed sliding motion
relating to a linear switching manifold
in the case that system is the same in directions tangent to the switching manifold.
In this case, the equation for motion
in the direction orthogonal to the switching manifold, $x$,
is independent of the variable representing displacement tangent to the
switching manifold, $y$,
and for this reason is amenable to an exact analysis.
We leave a description of more general scenarios for subsequent work.
In \S\ref{sec:X} we analyze the stochastic differential equation for $x(t)$
and derive its quasi-steady-state distribution.
In \S\ref{sec:Y} we analyze the equation for $y(t)$
and obtain expressions for the mean and variance of $y(t)$.
Derivations for this section are given in \S\ref{sec:CALC} and Appendix \ref{sec:PROOFS}.
Conclusions are presented in \S\ref{sec:CONC}.

\section{Periodic orbits with sliding and relay control}
\label{sec:MOTIV}

Periodic orbits involving sliding have recently been described in models
of relay control systems.
Broadly speaking, a relay control system is a system that aims to control a variable
using the measurements of an input signal via a switching action
\cite{Ts84,FrPo02,DoBi01,AsMu08}\removableFootnote{
One of the most fundamental problems in control theory
is that of determining the stability of the solution.
For this purpose multiple Lyapunov functions have recently been developed 
for the case of switched or hybrid control systems \cite{Br98,Li03}.
For systems that incorporate noise or additional functions representing uncertainty,
Lyapunov functions have been used to determine
generalized notions of stability \cite{McGl90,Zh98}.
In addition the robustness of the solution to noise is a desired property.
There seem to be only isolated results for the stability and robustness
of switched or hybrid control systems incorporating
a Wiener process \cite{RaMi10,FeZh06,ChLi04}
or an arbitrary uncertainty function \cite{SkEv99}.
However the listed results apply only to equilibrium solutions
whereas in regards to this paper we are more interested in the stability
and robustness of periodic solutions.
Also some recent work has been done to develop robust switching control laws
\cite{MhEl05,HuXu99,Su04}.
}.
Relay control systems are commonly modelled by
\begin{equation}
\begin{split}
\dot{\bx} &= A\bx + Bu \;, \\
\varphi &= C^{\sf T} \bx \;, \\
u &= -{\rm sgn}(\varphi) \;,
\end{split}
\label{eq:relayControlSystem}
\end{equation}
where $\bx \in \mathbb{R}^N$,
$\varphi$ is the signal measurement
and $u$ is the control response, \cite{Jo03,ZhMo03,DiBu08}.
The system (\ref{eq:relayControlSystem}) is a Filippov system
with a single switching manifold,
$\{ \bx ~|~ C^{\sf T} \bx = 0 \}$,
on which sliding may occur.
In this system sliding corresponds to the idealized scenario of
discrete switching events occurring continuously in time.
Periodic orbits of (\ref{eq:relayControlSystem}) that involve sliding are described in
\cite{DiJo01,JoRa99,JoBa02,ZhFe10}.

\begin{figure}[t!]
\begin{center}
\setlength{\unitlength}{1cm}
\begin{picture}(16,13.5)
\put(2.75,6.5){\includegraphics[height=7cm]{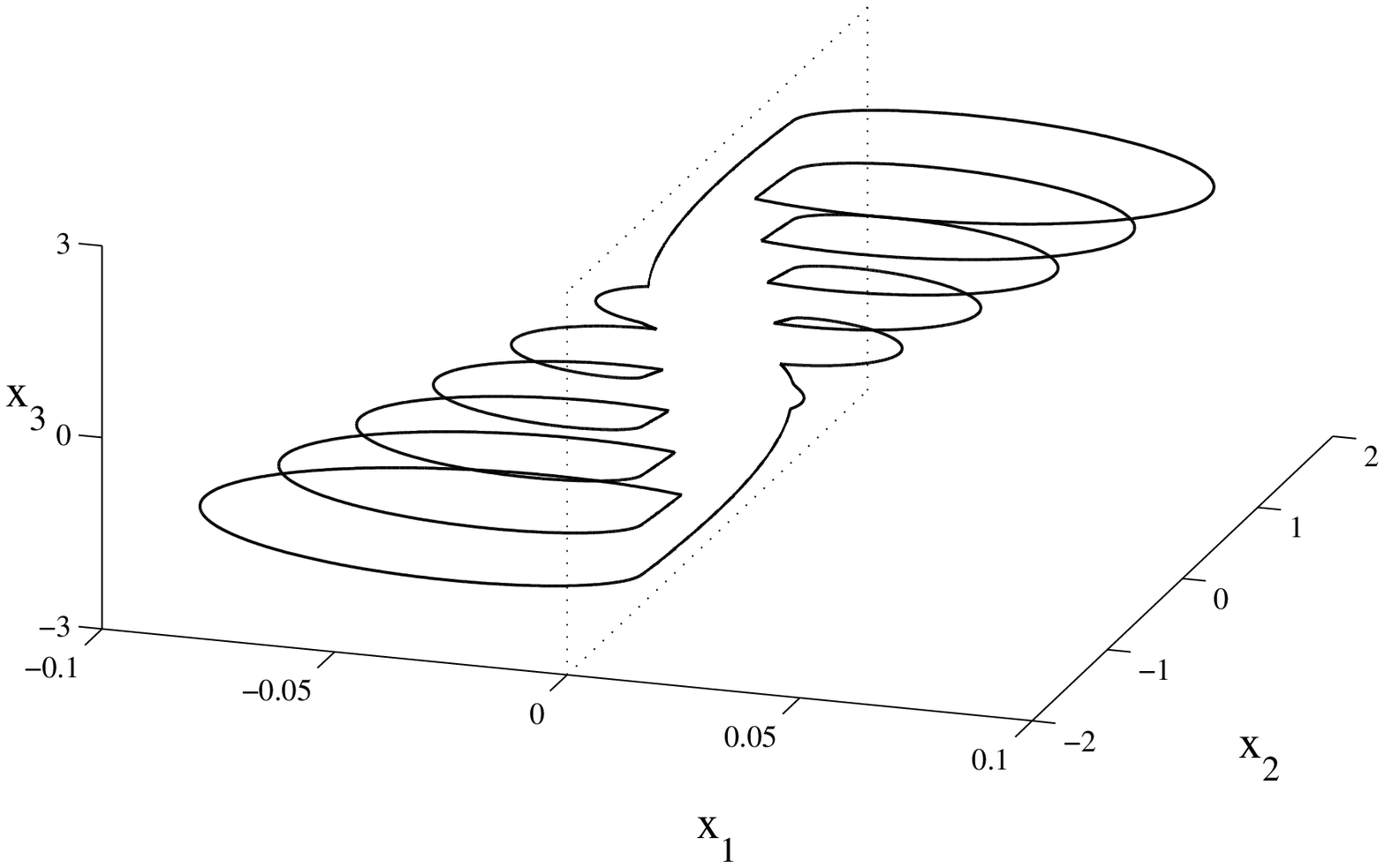}}	
\put(0,0){\includegraphics[height=6cm]{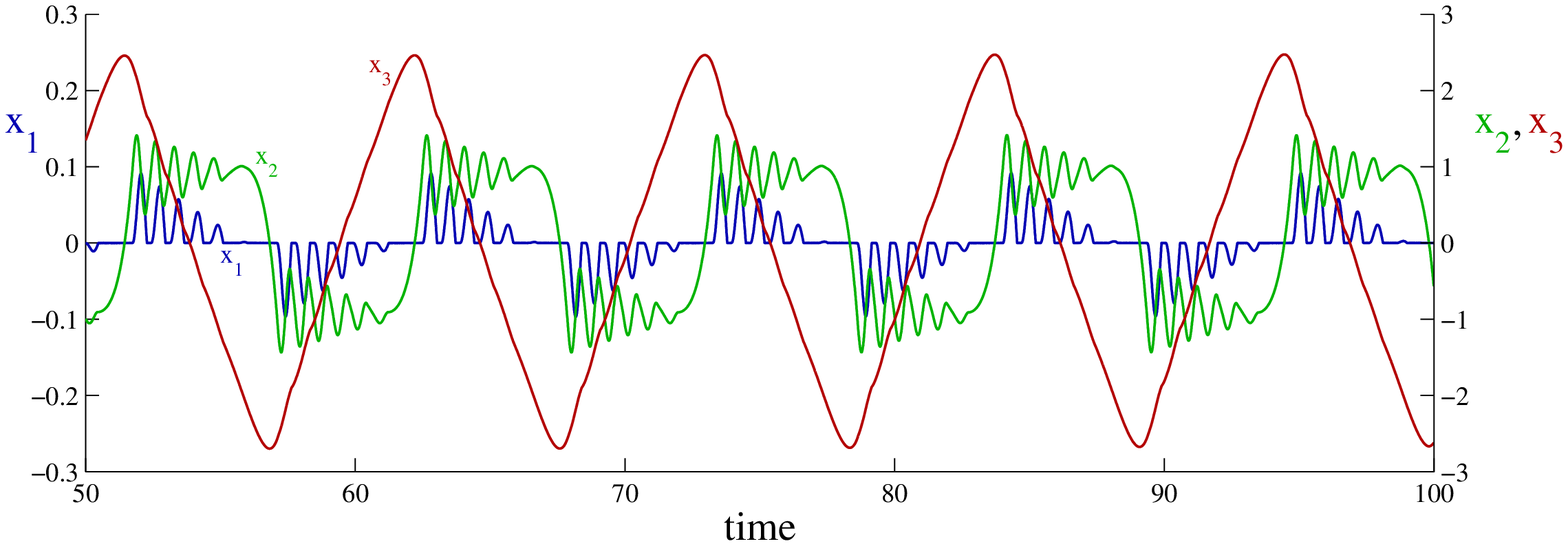}}		
\put(3.5,12){\large \sf \bfseries A}
\put(1.5,6){\large \sf \bfseries B}
\end{picture}
\caption{
A stable periodic orbit of (\ref{eq:relayControlSystem}) with
(\ref{eq:ABCvalues}) and $\zeta = -0.06$.
Here $\bx = (x_1, x_2, x_3)^{\sf T}$.
The periodic orbit exhibits sliding on the switching manifold, $x_1 = 0$.
\label{fig:noisyExampleDet}
}
\end{center}
\end{figure}

As an example we consider the following canonical form,
taken from \cite{DiJo01} (also given in \cite{DiBu08}),
\begin{equation}
A = \left[ \begin{array}{ccc}
-20 \zeta - \frac{1}{20} & 1 & 0 \\
-\zeta - 100 & 0 & 1 \\
-5 & 0 & 0
\end{array} \right] \;, \qquad
B = \left[ \begin{array}{c} 1 \\ -2 \\ 1 \end{array} \right] \;, \qquad
C = \left[ \begin{array}{c} 1 \\ 0 \\ 0 \end{array} \right] \;,
\label{eq:ABCvalues}
\end{equation}
where $\zeta \in \mathbb{R}$ is a parameter\removableFootnote{
The matrix $A$ may be assumed to be a companion matrix.
Three free parameters suffice to describe all possible $3 \times 3$ companion matrices.
The particular form
\begin{equation}
A = \left[ \begin{array}{ccc}
-2 \zeta \omega - \lambda & 1 & 0 \\
-2 \zeta \omega \lambda - \omega^2 & 0 & 1 \\
-\lambda \omega^2 & 0 & 0
\end{array} \right] \;,
\label{eq:Acompanion}
\end{equation}
has been used because it yields a conveniently simple transfer function.
The eigenvalues of (\ref{eq:Acompanion}) are
$-\lambda$ and $-\omega \zeta \pm {\rm i} |\omega| \sqrt{1 - \zeta^2}$.
Oddly both \cite{DiJo01} and \cite{DiBu08} state the complex eigenvalues incorrectly
(pg.~1134 of \cite{DiJo01} gives $-\omega \zeta \pm {\rm i} \omega \sqrt{1 - \zeta}$,
and pg.~364 of \cite{DiBu08} gives $\zeta \pm {\rm i} \omega$).
}.
Fig.~\ref{fig:noisyExampleDet} illustrates a stable periodic orbit of
(\ref{eq:relayControlSystem}) with (\ref{eq:ABCvalues})
involving sliding motion on the switching manifold.
There are 12 separate sliding segments per period.
These correspond to time intervals for which $x_1$
(the first component of the vector, $\bx$) is zero.
The stability of periodic orbits with sliding may be determined
by analyzing the Jacobian of a return map \cite{DiJo01,ZhFe10,TaOs09}.
The robustness of periodic orbits with sliding
has been briefly investigated by studying the size of the basin of
attraction of the periodic orbit \cite{TaOs09},
and imposing a short time between consecutive switching events \cite{DiJo02}.

\begin{figure}[t!]
\begin{center}
\setlength{\unitlength}{1cm}
\begin{picture}(16,12.7)
\put(0,6.7){\includegraphics[height=6cm]{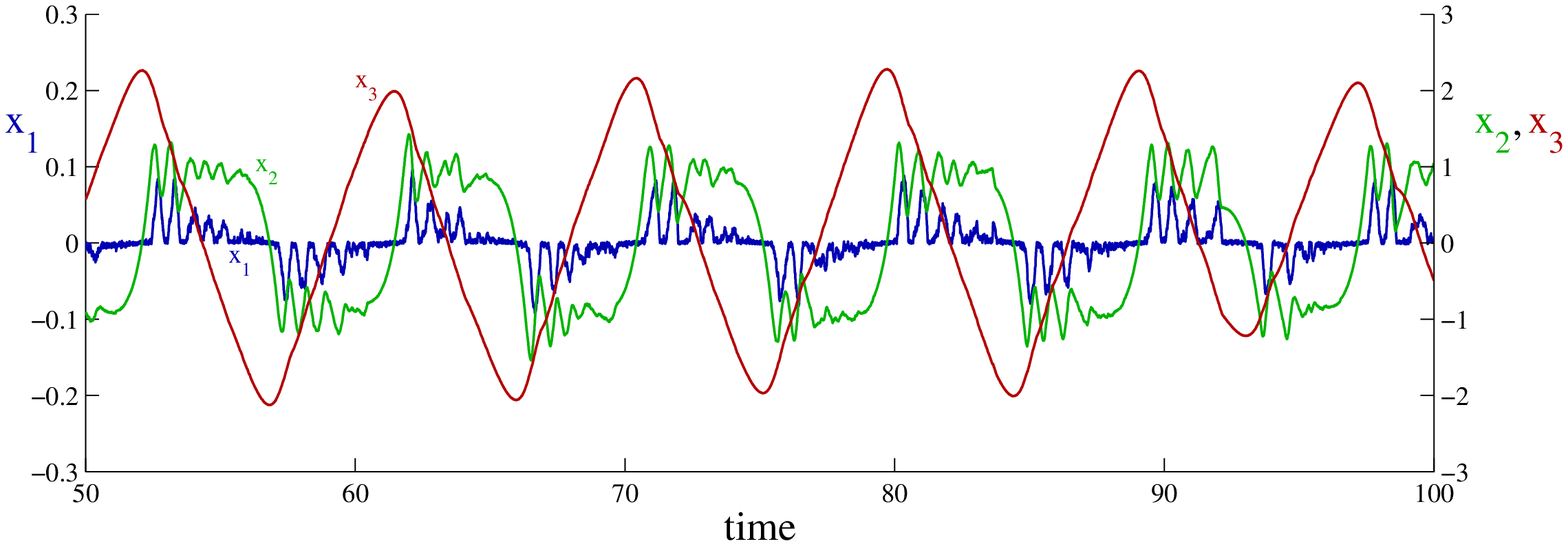}}
\put(0,.2){\includegraphics[height=6cm]{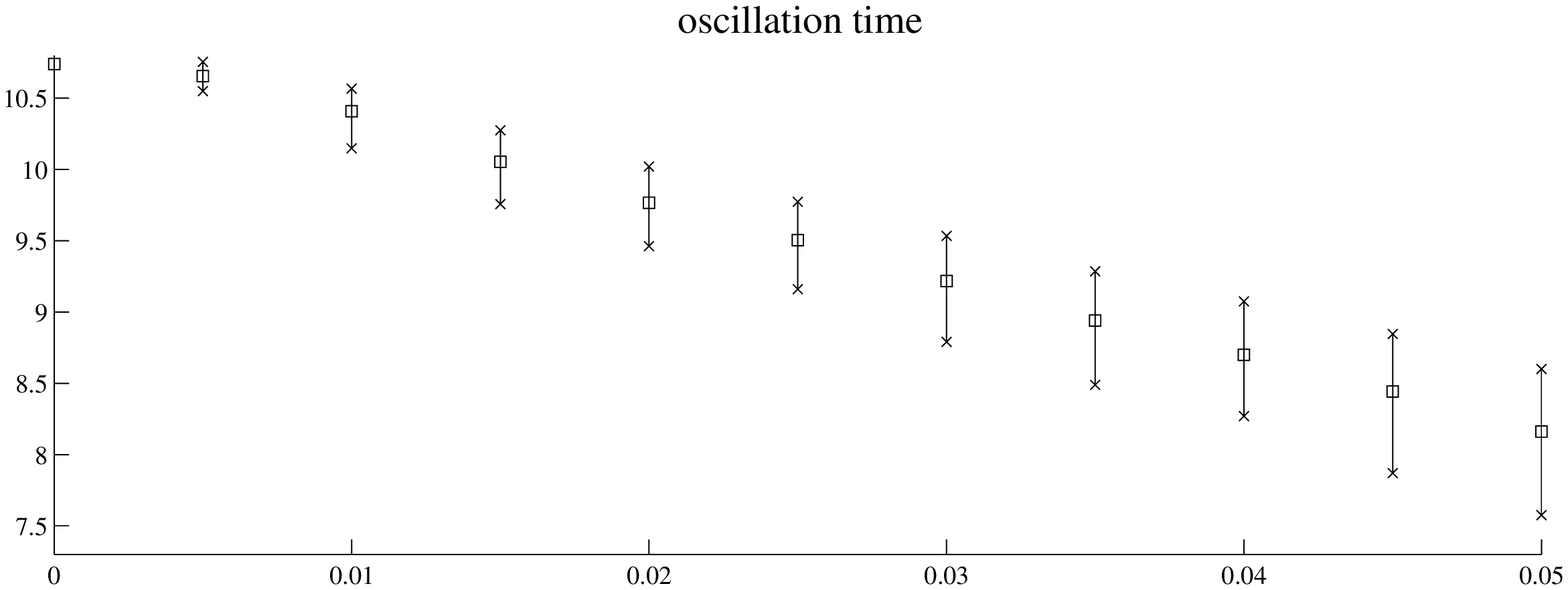}}
\put(8,0){$\sqrt{\ee}$}
\put(1.4,12.2){\large \sf \bfseries A}
\put(1,5.7){\large \sf \bfseries B}
\end{picture}
\caption{
Panel A shows a time series of (\ref{eq:relayControlSystem3}) with $\ee = 0.001$,
using the same parameter values as in Fig.~\ref{fig:noisyExampleDet}.
Panel B shows the median, upper quartile and lower quartile values of
1000 numerically computed oscillation times, for several values of $\ee$.
To obtain each oscillation time,
we computed an orbit up to $t = 100$ and identified the last
three instances at which the value of $x_3$ changed sign
(discounting rapid sign changes over a handful of grid points near this value due to noise),
then subtracted the first time from the third time.
Orbits were computed with the
Euler-Maruyama method of fixed step size, $\Delta t = 0.0001$.
\label{fig:noisyExample}
}
\end{center}
\end{figure}

Fig.~\ref{fig:noisyExample}-A
shows a typical orbit for the system when noise
is added to the control signal as
\begin{equation}
d\bx = \left( A \bx - B \,{\rm sgn} \left( C^{\sf T} \bx \right) \right) \, dt +
\sqrt{\ee} B \, dW(t) \;.
\label{eq:relayControlSystem3}
\end{equation}
By comparing with Fig.~\ref{fig:noisyExampleDet}-B,
this shows that the noise may dampen the oscillations and induce
an increase in frequency\removableFootnote{
In future I could use a more sophisticated numerical method
such as the split-step backward method described in \cite{BaTa12}.
}.
Fig.~\ref{fig:noisyExample}-B shows that the
average length of time per oscillation of $x_3$
decreases as the noise amplitude increases.
We have observed similar behaviour for different values of $\zeta$.
In \S\ref{sec:CONC} we use the results below
to speculate on the reason for this behaviour.

\section{A simple set of equations for stochastically perturbed sliding motion}
\label{sec:TOY}

We are interested in the effects of noise on the dynamics of (\ref{eq:Filippov})
near a switching manifold.
In this section we introduce a 
simple set of equations that approximates (\ref{eq:Filippov}) near a switching manifold,
then formulate the inclusion of noise.

If a switching manifold of (\ref{eq:Filippov})
is locally $C^K$ we can choose a coordinate system
such that in a neighbourhood of the origin, $\bx = 0$,
the switching manifold is simply $e_1^{\sf T} \bx = O(K)$ \cite{DiBu01}.
In this paper we are not concerned with effects
due to nonsmoothness in the switching manifold
and for this reason suppose that the switching manifold is
the coordinate plane, $e_1^{\sf T} \bx = 0$.
Two smooth subsystems govern the nearby flow,
thus, locally, we may write the deterministic system as
\begin{equation}
\dot{\bx} = \left\{ \begin{array}{lc}
F^{(L)}(\bx) \;, & e_1^{\sf T} \bx < 0 \\
F^{(R)}(\bx) \;, & e_1^{\sf T} \bx > 0
\end{array} \right. \;,
\label{eq:Filippov2}
\end{equation}
where $F^{(L)}$ and $F^{(R)}$ are, say, $C^1$.

Suppose there exists a section of the switching manifold,
call it $\Sigma$, for which
$e_1^{\sf T} F^{(L)}(\bx) > 0$ and $e_1^{\sf T} F^{(R)}(\bx) < 0$,
as in Fig.~\ref{fig:noisyFlowExSlide}-A.
In this scenario, forward orbits arrive at $\Sigma$ from either side.
We use Filippov's definition
to define dynamics constrained to $\Sigma$ \cite{Fi60,Fi88,DiBu08}.
$\Sigma$ is known as an {\em attracting sliding region} and
evolution on $\Sigma$ is referred to as {\em sliding motion}.

However, with additive noise, the system (\ref{eq:Filippov2})
seems to be too complex for us analyze to a degree of detail that is useful.
This is because both $F^{(L)}$ and $F^{(R)}$
depend on all components of the vector $\bx$,
and we have been unable to analytically solve the resulting $N$-dimensional
stochastic differential equation.
Consequently, for this paper,
which represents a first detailed analysis of sliding motion with noise,
we ignore the dependency of $F^{(L)}$ and $F^{(R)}$
on components of $\bx$ parallel to $\Sigma$
as this enables us to reduce mathematical problems to one dimension
but still capture what seems to be the essence of 
stochastically perturbed sliding motion.
Moreover, this provides a useful approximation to the general case
over short time-frames.

Given that $F^{(L)}$ and $F^{(R)}$ are functions of only $e_1^{\sf T} \bx$,
the remaining $N-1$ components of $\bx$ may be treated identically and for this
reason it suffices to study a two-dimensional system.
We let $\bx = [x~y]^{\sf T}$ and add
small amplitude, white, Gaussian noise
independent to the state of the system.
Assuming for simplicity that the noise in $x$ is independent of the noise in $y$,
the resulting stochastic differential equation may be written as
\begin{eqnarray}
\left[ \begin{array}{c} dx \\ dy \end{array} \right]
&=& \left[ \begin{array}{c} \phi(x) \\ \psi(x) \end{array} \right] \, dt
+ \sqrt{\ee} \,
\left[ \begin{array}{c} dW_1(t) \\ \sqrt{\kappa} \, dW_2(t) \end{array} \right] \;,
\label{eq:e}
\end{eqnarray}
where $W_1(t)$ and $W_2(t)$ are independent Brownian motions,
$0 < \ee \ll 1$ and $\kappa > 0$ are constants,
and $\phi$ and $\psi$ are piecewise-$C^1$ that
for small $|x|$ are given by
\begin{eqnarray}
\phi(x) &=& \left\{ \begin{array}{lc}
a_L + c_L x + o(|x|) \;, & x<0 \\
-a_R + c_R x + o(|x|) \;, & x>0
\end{array} \right. \;,
\label{eq:phi} \\
\psi(x) &=& \left\{ \begin{array}{lc}
b_L + d_L x + o(|x|) \;, & x<0 \\
b_R + d_R x + o(|x|) \;, & x>0
\end{array} \right. \;.
\label{eq:psi}
\end{eqnarray}
We assume
\begin{equation}
a_L, a_R > 0 \;,
\label{eq:aLaRpos}
\end{equation}
to ensure that in the absence of noise
the switching manifold ($x=0$) is an attracting sliding region.
Since $\phi$ and $\psi$ are independent of $\ee$,
their coefficients are $O(1)$.
Consequently, for $x(0)$ near zero, orbits of (\ref{eq:e})
likely remain near $x=0$ for relatively long periods of time,
as shown in \S\ref{sub:EXIT},
and for this reason we
do not specify the behaviour of $\phi$ and $\psi$ for large $|x|$.

\section{Properties of $x(t)$}
\label{sec:X}

Since (\ref{eq:e}) lacks dependency on $y$,
the equation for $dx$ is decoupled from $y$:
\begin{eqnarray}
dx
&=& \phi(x) \, dt + \sqrt{\ee} \, dW_1(t) \nonumber \\
&=& \left\{ \begin{array}{lc}
a_L + c_L x + o(|x|) \;, & x<0 \\
-a_R + c_R x + o(|x|) \;, & x>0
\end{array} \right\} \, dt
+ \sqrt{\ee} \, dW_1(t) \;.
\label{eq:dx}
\end{eqnarray}
Given $x(0) = x_0$,
let $p_\ee(x,t|x_0)$ denote the transitional \pdf~for
the value of $x(t)$, as governed by (\ref{eq:dx}).
Despite the discontinuity at $x=0$,
$p_\ee(x,t|x_0)$ is unique and continuous
on $\mathbb{R} \times (0,\infty) \times \mathbb{R}$.
For $x \ne 0$ and $t > 0$, the \pdf~satisfies the Fokker-Planck equation
\begin{equation}
\frac{\partial p_\ee}{\partial t} =
- \frac{\partial (\phi p_\ee)}{\partial x} 
+ \frac{\ee}{2} \frac{\partial^2 p_\ee}{\partial x^2} \;,
\label{eq:fp}
\end{equation}
with the initial condition $p_\ee(x,0|x_0) = \delta(x-x_0)$, \cite{Sc80,Ga85,Ri84}.
In \S\ref{sec:CALC} we provide an explicit expression for
$p_\ee(x,t|x_0)$ in the special case that $\phi$ is piecewise-constant.
For the remainder of this section we use (\ref{eq:fp})
to determine the long time behaviour of $x(t)$.

If $\phi(x) > 0$ for all $x < 0$, and $\phi(x) < 0$ for all $x > 0$,
then (\ref{eq:dx}) has a steady-state density on $\mathbb{R}$
centered about the origin.
Otherwise, $\phi(x) = 0$ for some $x \ne 0$, and
with nonzero probability
orbits may cross this value of $x$ and
undergo dynamics far from the origin not described
by the expansion (\ref{eq:phi}).
However, regardless of the global nature of $\phi$,
since $\ee$ is small the local attraction to the origin is relatively strong.
Thus we expect orbits to remain near the origin for long periods of time
and be distributed by a quasi-steady-state distribution for large but finite $t$.
In \S\ref{sub:EXIT} we determine the mean escape time
of orbits from an $O(1)$ neighbourhood of the origin.
In \S\ref{sub:QSS} we use (\ref{eq:fp}) to derive
the quasi-steady-state probability density function asymptotically.

\subsection{Escape from a neighbourhood of $x=0$}
\label{sub:EXIT}

For the function $\phi(x)$, (\ref{eq:phi}), with (\ref{eq:aLaRpos}),
in the case that $\phi(x) = 0$ for some $x \ne 0$,
it is necessary to identify a value $x_b > 0$, independent of $\ee$, such that
\begin{equation}
\min_{|x| \le x_b} |\phi(x)| \ge \frac{1}{2} \min (a_L,a_R) \;.
\end{equation}
Then in the interval $[-x_b,x_b]$, the drift of (\ref{eq:dx}) is towards $x=0$.
With $\ee > 0$ and any initial condition $x_0 \in (-x_b,x_b)$,
$x(t)$ will eventually escape $[-x_b,x_b]$ with probability 1.
Calculating the time to escape is a standard problem
in the context of a single potential well, 
where the potential function is given by 
\begin{equation}
U(x) = -\int_0^x \phi(y) \,dy \;.
\end{equation}
The mean escape time, $\overline{T}(x_0)$, may be found
exactly \cite{Sc10,GrVa99,Ga85}\removableFootnote{
DJWS: Specifically, $\overline{T}(x_0)$ satisfies Dynkin's equation:
\begin{equation}
\frac{\ee}{2} \overline{T}'' - U'(x_0) \overline{T}' = -1 \;,
\end{equation}
on $[-x_b,x_b]$, and has homogeneous boundary conditions, $\overline{T}(\pm x_b) = 0$.
By using an integrating factor we arrive at the following integral expression:
\begin{eqnarray}
\overline{T}(x_0) &=& \frac{2}{\ee} \int_{-x_b}^{x_0} {\rm e}^{\frac{2 U(z)}{\ee}}
\left( C - \int_{-x_b}^z {\rm e}^{-\frac{2 U(y)}{\ee}} \,dy \right) \,dz \;, \\
{\rm where~~} C &=&
\frac{\int_{-x_b}^{x_b} {\rm e}^{\frac{2 U(\hat{z})}{\ee}}
\int_{-x_b}^{\hat{z}} {\rm e}^{-\frac{2 U(y)}{\ee}} \,dy \,d\hat{z}}
{\int_{-x_b}^{x_b} {\rm e}^{\frac{2 U(\hat{z})}{\ee}} \,d\hat{z}} \;,
\end{eqnarray}
see Eqn.~6.54 of \cite{Sc10}, Eqn.~3.13 of \cite{GrVa99} and Eqn.~5.2.158 of \cite{Ga85}.
}.
Via Laplace's method of asymptotic evaluation of integrals \cite{BeOr78},
it follows that whenever $|x_0| < x_b$\removableFootnote{
We do not need to specify 
a closed subinterval of $[-x_b,x_b]$ for the domain of $x_0$
because our statements are asymptotic in $\ee$.
Since the statement $|x_0| < x_b$ is independent of $\ee$,
it implies that $x_b - |x_0| = O(1) \gg \ee > 0$.
The point is that we are essentially choosing $\ee$ after we have chosen $x_0$.
}
there exist $\ee$-independent constants $\alpha_1$ and $\alpha_2$ such that\removableFootnote{
Ideally we would like to have a precise asymptotic
expression for the probability that an orbit has escaped, $G(t,x_0)$,
but this is probably more effort than is necessary for the purposes of this paper.
A standard result in this direction is that if we write
$G(t,x_0) = \sum_{n=1}^{\infty} \beta_n {\rm e}^{-\lambda_n t} X_n(x_0)$,
as obtained by using separation of variables to solve the
backwards Fokker-Planck equation,
then the eigenvalues satisfy $0 < \lambda_1 < \lambda_2 < \cdots$,
and $\lambda_1 \sim \frac{1}{\overline{T}(x_0)}$, \cite{Sc10}.
An approximation to the first eigenfunction, $X_1$,
is described on pg.~352 of \cite{Sc10}.
}
\begin{equation}
\overline{T}(x_0) \sim \ee \alpha_1 {\rm e}^{\frac{\alpha_2}{\ee}} \;.
\end{equation}
For instance if $U(x_b) < U(-x_b)$ (other cases are similar),\removableFootnote{
In this case escape is far is more likely to occur through $x_b$.
Note that $a_L$ and $a_R$ arise from an asymptotic evaluation of the constant, $C$.
(In the special case that $\phi$ is piecewise-constant,
I derived an exact non-integral expression for $\overline{T}(x_0)$
and found it to be consistent with this equation.)
By symmetry, the result in the case $U(x_b) > U(-x_b)$ is analogous;
the result for the case $U(x_b) = U(-x_b)$ is also similar,
but the coefficients are messier.
}
\begin{equation}
\overline{T}(x_0) \sim \frac{\ee (a_L+a_R)}{2 a_L a_R U'(x_b)} {\rm e}^{\frac{2 U(x_b)}{\ee}} \;.
\end{equation}

\subsection{The quasi-steady-state probability density function}
\label{sub:QSS}

In \S\ref{sub:EXIT} we showed that the mean
escape time from an $\ee$-independent neighbourhood $[-x_b,x_b]$
is exponentially large in $\frac{1}{\ee}$.
Consequently, we can assume that the probability an orbit escapes $[-x_b,x_b]$
within the polynomial time $\frac{1}{\ee^{M}}$, for any fixed $M > 0$,
is extremely small.
We let
\begin{equation}
\check{t} = \ee^M t \;,
\label{eq:tcheck}
\end{equation}
represent the long time scale, and look for a solution
to the Fokker-Planck equation (\ref{eq:fp}) as a function of $x$ and $\check{t}$.
By substituting (\ref{eq:tcheck}) with the
WKB-type expansion \cite{Sc80,Ga85},
\begin{equation}
p_\ee(x,t|x_0) = {\rm e}^{\frac{q_\ee(x,t|x_0)}{\ee}} \;,
\label{eq:wkb}
\end{equation}
into (\ref{eq:fp}), we arrive at
\begin{equation}
\ee^M \frac{\partial q_\ee}{\partial \check{t}} =
\left( \frac{\partial q_\ee}{\partial x} + \ee \frac{\partial}{\partial x} \right)
\left( \frac{1}{2} \frac{\partial q_\ee}{\partial x} - \phi \right) \;.
\end{equation}
Therefore
\begin{equation}
\frac{\partial q_\ee(x,\ee^{-M} \check{t}|x_0)}{\partial x}
= 2 \phi(x) + O(\ee^M) \;,
\end{equation}
and so by integrating $\phi(x)$ we obtain
\begin{equation}
q_\ee(x,\ee^{-M} \check{t}|x_0) =
r(\check{t}) + \left\{ \begin{array}{lc}
2 a_L x + c_L x^2 + o(x^2) \;, & x \le 0 \\
-2 a_R x + c_R x^2 + o(x^2) \;, & x \ge 0
\end{array} \right\} + O(\ee^M) \;,
\label{eq:q}
\end{equation}
for an $\ee$-independent function $r$.
As a function of $x$ and $t$, the dependence of this solution
on $t$ and $x_0$ appears only in $O(\ee^M)$ terms
which may be ignored.
Consequently we treat the solution as solely a function of $x$
and refer to it as the {\em quasi-steady-state solution}, $p_{{\rm qss},\ee}(x)$.
Specifically (\ref{eq:wkb}) and (\ref{eq:q}) combine to give
\begin{equation}
p_{{\rm qss},\ee}(x) = \left\{ \begin{array}{lc}
\frac{K_\ee}{\ee} {\rm e}^{\frac{1}{\ee}
( 2 a_L x + c_L x^2 + o(x^2) ) + O(\ee^M)} \;, & x \le 0 \\
\frac{K_\ee}{\ee} {\rm e}^{\frac{1}{\ee}
( -2 a_R x + c_R x^2 + o(x^2) ) + O(\ee^M)} \;, & x \ge 0
\end{array} \right. \;,
\label{eq:qss}
\end{equation}
where we must have\removableFootnote{
I arrived at this by differentiating the integral of
$p_{{\rm qss},\ee}$ over $[-\sqrt{\ee},\sqrt{\ee}]$
and with respect to $\ee$ to obtain
the Taylor series of $K_\ee$ centred at $\ee = 0$.
}
\begin{equation}
K_\ee = \frac{2 a_L a_R}{a_L+a_R} -
\frac{a_L^3 c_R + a_R^3 c_L}{a_L a_R (a_L+a_R)^2} \ee + o(\ee) \;,
\label{eq:K}
\end{equation}
to ensure $p_{{\rm qss},\ee}$ is normalized.

For small $\ee$ and $x_0$, the transitional \pdf~of (\ref{eq:dx}),
$p_\ee(x,t|x_0)$, quickly settles to (\ref{eq:qss}).
The scaling
\begin{equation}
\hat{x} = \frac{x}{\ee} \;, \qquad
\hat{t} = \frac{t}{\ee} \;,
\label{eq:xthat}
\end{equation}
transforms (\ref{eq:dx}) to
\begin{equation}
d\hat{x} = \phi(\ee \hat{x}) \,d\hat{t} + dW_1(\hat{t}) \;,
\label{eq:xtScaling}
\end{equation}
from which we infer that $p_\ee(x,t|x_0)$ approaches (\ref{eq:qss})
on an $O(\ee)$ time-scale, when $x_0 = O(\ee)$.
Furthermore, for times in the range $\ee^{1-\delta} \le t \le \ee^{-M}$,
where $\delta > 0$,
it is reasonable to suppose $x \sim p_{{\rm qss},\ee}$,
in which case
\begin{eqnarray}
\big< {\rm sgn}(x) \big> &=&
\frac{a_L-a_R}{a_L+a_R} +
\frac{a_L^2 c_R - a_R^2 c_L}{a_L a_R (a_L+a_R)^2} \ee + o(\ee) \;, 
\label{eq:sgnxMean} \\
\big< x \big> &=&
\frac{a_L-a_R}{2 a_L a_R} \ee + O(\ee^2) \;,
\label{eq:xMean} \\
\big< x \, {\rm sgn}(x) \big> &=&
\frac{a_L^2 + a_R^2}{2 a_L a_R (a_L+a_R)} \ee + O(\ee^2) \;,
\label{eq:xsgnxMean}
\end{eqnarray}
which are useful in the next section.

\section{Moments of $y(t)$}
\label{sec:Y}

In this section we compute the mean of $y(t)$
and conjecture the leading order term of its variance.
We assume $x \sim p_{{\rm qss},\ee}$ at all times under consideration
which greatly simplifies calculations.
We begin by deriving $y(t)$ when $\ee = 0$.

\subsection{Deterministic sliding motion}

When $\ee = 0$, (\ref{eq:e}) is the 
Filippov system:
\begin{equation}
\left[ \begin{array}{c} \dot{x} \\ \dot{y} \end{array} \right] =
\left[ \begin{array}{c} \phi(x) \\ \psi(x) \end{array} \right] \;.
\label{eq:ode}
\end{equation}
As in \cite{Fi60,Fi88,DiBu08},
we define a vector field for sliding motion
on the switching manifold ($x=0$)
by the unique convex combination of the two vector fields at the manifold
that is tangent to the manifold.
That is,
\begin{equation}
\left[ \begin{array}{c} \dot{x} \\ \dot{y} \end{array} \right]_{\rm slide} =
(1-q) \left[ \begin{array}{c} a_L \\ b_L \end{array} \right] +
q \left[ \begin{array}{c} -a_R \\ b_R \end{array} \right] \;,
\end{equation}
for the unique scalar, $q \in (0,1)$,
for which $\dot{x}_{\rm slide} = 0$.
Solving $\dot{x}_{\rm slide} = 0$ gives $q = \frac{a_L}{a_L+a_R}$ and therefore
\begin{equation}
\dot{y}_{\rm slide} = \frac{a_R b_L + a_L b_R}{a_L+a_R} \;.
\end{equation}
Consequently, if $(x(0),y(0)) = (0,y_0)$,
then $x_{\rm slide}(t) \equiv 0$ and
\begin{equation}
y_{\rm slide}(t) \equiv y_0 + \frac{a_R b_L + a_L b_R}{a_L+a_R} t \;.
\label{eq:ySlide}
\end{equation}

\subsection{The mean of $y(t)$}

From (\ref{eq:e}) and (\ref{eq:psi}) we have
\begin{eqnarray}
dy &=& \left( \frac{b_L+b_R}{2} - \frac{b_L-b_R}{2} \, {\rm sgn}(x(t))
+ \frac{d_L+d_R}{2} x(t)
- \frac{d_L-d_R}{2} x(t) \, {\rm sgn}(x(t)) + o(|x(t)|) \right) \, dt \nonumber \\
&&+~\sqrt{\ee \kappa} \, dW_2(t) \;.
\label{eq:dy}
\end{eqnarray}
Integration yields
\begin{eqnarray}
y(t) &=& y_0 + \frac{b_L+b_R}{2} t
- \frac{b_L-b_R}{2} \int_0^t {\rm sgn}(x(s)) \, ds
+ \frac{d_L+d_R}{2} \int_0^t x(s) \, ds \\ \nonumber
&&-~\frac{d_L-d_R}{2} \int_0^t x(s)
\, {\rm sgn}(x(s)) \, ds
+ \int_0^t o(|x(s)|) \, ds
+ \sqrt{\ee \kappa} \, W_2(t) \;,
\label{eq:y}
\end{eqnarray}
and therefore
\begin{eqnarray}
\big< y(t) \big> &=& y_0 + \frac{b_L+b_R}{2} t
- \frac{b_L-b_R}{2} \int_0^t \big< {\rm sgn}(x(s)) \big> \, ds
+ \frac{d_L+d_R}{2} \int_0^t \big< x(s) \big> \, ds \\ \nonumber
&&-~\frac{d_L-d_R}{2} \int_0^t
\big< x(s) \, {\rm sgn}(x(s)) \big> \, ds
+ \int_0^t \big< o(|x(s)|) \big> \, ds \;.
\label{eq:yMeanGen}
\end{eqnarray}
If $x \sim p_{{\rm qss},\ee}$,
by substituting (\ref{eq:sgnxMean})-(\ref{eq:xsgnxMean}) we obtain
\begin{equation}
\big< y(t) \big> = y_{\rm slide}(t) +
\frac{(a_L^2 d_R - a_R^2 d_L)(a_L+a_R) -
(a_L^2 c_R - a_R^2 c_L)(b_L-b_R)}
{2 a_L a_R (a_L+a_R)^2} \ee t + o(\ee) \;,
\label{eq:yMean}
\end{equation}
where the $\ee$-independent terms
have combined to form $y_{\rm slide}(t)$, (\ref{eq:ySlide}).

Therefore as $\ee \to 0$,
the mean of $y(t)$ limits on Filippov's sliding solution, $y_{\rm slide}(t)$.
This is non-trivial because Filippov's method, to obtain (\ref{eq:ySlide}),
and standard stochastic dynamical systems definitions, to obtain (\ref{eq:yMean}),
are not immediately related\removableFootnote{
It seems clear that in general
as $\ee \to 0$ the stochastic solution limits on Filippov's solution,
with reasonable assumptions.
However I am not aware of any other work that states this fundamental result.
My inclination is that a proof would require a technical measure-theoretic
approach and that it is likely to have already been proved, but just
locked up within some general theorems
on the existence, uniqueness and continuity of solutions to stochastic differential equations.
Filippov \cite{Fi60} shows that his sliding solution of a discontinuous vector field
equals the solution to a limiting mollified vector field.
The same is true, with certain assumptions,
for stochastic differential equations with discontinuous drift \cite{Co71}.
What I want is the case that the diffusion coefficient goes to zero in the limit.
This is done for smooth systems in, for instance, \cite{FrWe84}.
}.
Note that the perturbation of $\big< y(t) \big>$ from $y_{\rm slide}(t)$ is order $\ee$,
mirroring the result for smooth systems, see \S\ref{sec:INTRO}.
The explicit expression for the coefficient of the $O(\ee)$-term in
(\ref{eq:yMean}) is particularly useful.
For instance, we can see that if $c_L = c_R$ and $d_L = d_R$,
then we require the asymmetry $a_L \ne a_R$ in order for the $O(\ee)$-term
to be nonzero.

\subsection{The variance of $y(t)$}

\begin{figure}[t!]
\begin{center}
\setlength{\unitlength}{1cm}
\begin{picture}(10,5.5)
\put(0,0){\includegraphics[height=5cm]{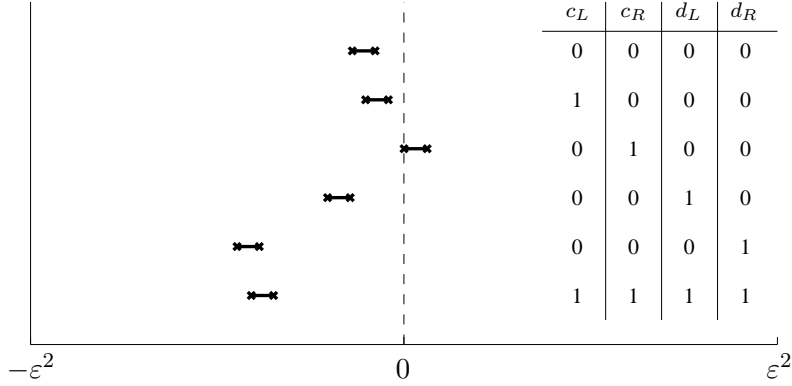}}
\put(-.25,0){\small $-\ee^2$}
\put(4.92,0){\small $0$}
\put(9.85,0){\small $\ee^2$}
\put(7.18,4.8){\scriptsize $c_L$}
\put(7.91,4.8){\scriptsize $c_R$}
\put(8.64,4.8){\scriptsize $d_L$}
\put(9.37,4.8){\scriptsize $d_R$}
\end{picture}
\caption{
The six horizontal bars are $95\%$ confidence intervals for
${\rm Var}(y(t)) - \left( \kappa + \frac{(b_L-b_R)^2}{(a_L+a_R)^2} \right) \ee t$,
when $(a_L,a_R,b_L,b_R) = (2,1,1,0)$, $\kappa = 0$, $\ee = 0.01$, $t = 1$,
and $\phi$ and $\psi$ are piecewise-linear 
using different values of $c_L$, $c_R$, $d_L$ and $d_R$
as indicated by the table.
Each confidence interval was obtained from $10^6$
orbits computed numerically with the Euler-Maruyama method of
fixed step size, $\Delta t = 0.0001$.
The results are consistent with Conjecture \ref{cj:yVar}
which predicts that the values of $c_L$, $c_R$, $d_L$ and $d_R$
affect the magnitude of ${\rm Var}(y(t))$ by at most $O(\ee^2)$.
\label{fig:checkVarMany}
}
\end{center}
\end{figure}

The variance of $y(t)$ may be computed via\removableFootnote{
Specifically
\begin{eqnarray}
{\rm Var}(y(t)) &=& \ee \kappa t +
\frac{(b_L-b_R)^2}{4} \int_0^t \int_0^t
\big< {\rm sgn}(x(s)) \, {\rm sgn}(x(u)) \big> -
\big< {\rm sgn}(x(s)) \big>
\big< {\rm sgn}(x(u)) \big> \, ds \, du \nonumber \\
&&-~\frac{(b_L-b_R) (d_L+d_R)}{2} \int_0^t \int_0^t
\big< {\rm sgn}(x(s)) \, x(u) \big> -
\big< {\rm sgn}(x(s)) \big>
\big< x(u) \big> \, ds \, du \nonumber \\
&&+~\frac{(b_L-b_R) (d_L-d_R)}{2} \int_0^t \int_0^t
\big< {\rm sgn}(x(s)) x(u)
\, {\rm sgn}(x(u)) \big> - \big< {\rm sgn}(x(s)) \big>
\big< x(u) \, {\rm sgn}(x(u)) \big> \, ds \, du  \nonumber \\
&&+~\int_0^t \int_0^t
\big< {\rm sgn}(x(s)) o(|x(u)|) \big> -
\big< {\rm sgn}(x(s)) \big> \big< o(|x(u)|) \big> \, ds \, du \nonumber \\
&&+~\int_0^t \int_0^t
\big< O(|x(s)|) O(|x(u)|) \big> -
\big< O(|x(s)|) \big> \big< O(|x(u)|) \big> \, ds \, du \;.
\label{eq:yVar}
\end{eqnarray}
}
\begin{equation}
{\rm Var}(y(t)) = \big< y(t)^2 \big> - \big< y(t) \big>^2 \;,
\label{eq:varIdentity}
\end{equation}
however this requires knowledge of $p_\ee(x,t|x_0)$,
for which we have not been able to obtain a useful expression
in the case of general $\phi$.
We conjecture that the leading order terms of ${\rm Var}(y(t))$
are independent of non-constant terms in $\phi$ and $\psi$
because $x(t) = O(\ee)$ with high probability\removableFootnote{
In order to prove Conjecture \ref{cj:yVar}
there are two major technical hurdles to overcome,
neither of which I know how to deal with.
First, we must show that each double integral in (\ref{eq:yVar}),
except the first, is $O(\ee^2)$, or smaller.
Certainly, with the substitution $x = \ee \hat{x}$ (recall (\ref{eq:xthat}), (\ref{eq:xtScaling})),
from each double integral except the first,
we can factor out at least one $\ee$.
After the double integrals are simplified to single integrals,
a second $\ee$ comes out by
the corresponding time scaling: $u = \ee \hat{u}$.
Then, with $\ee^2$ successfully factored out,
the resulting double integrals are $O(1)$
because they now correspond to the ``hatted'' variables.
However I don't know how to show that the double integrals are finite
because this requires knowing that $\hat{p}_\ee \to \hat{p}_{{\rm qss},\ee}$
sufficiently quickly
(e.g.~show that
\begin{equation}
\int_0^{\frac{t}{\ee}} \big< {\rm sgn}(\hat{x}(0)) \hat{x}(\hat{u}) \big> -
\big< {\rm sgn}(\hat{x}(0)) \big> \big< \hat{x}(\hat{u}) \big> \, d\hat{u}
= \int_0^{\frac{t}{\ee}}
\int_{-\frac{x_b}{\ee}}^{\frac{x_b}{\ee}} 
{\rm sgn}(\hat{x}) p_{{\rm qss},\ee}(\hat{x})
\int_{\frac{-x_b}{\ee}}^{\frac{x_b}{\ee}}
\hat{y} \left( p_\ee(\hat{y},\hat{u}|\hat{x}) - p_{{\rm qss},\ee}(\hat{y}) \right)
\, d\hat{y} \, d\hat{x} \, d\hat{u}
\end{equation}
is finite).
Second, we must show that non-constant terms in the two pieces of $\phi$
contribute only higher order terms to the first double integral in (\ref{eq:yVar}).
}.
Indeed this is consistent with numerical simulations,
Fig.~\ref{fig:checkVarMany}\removableFootnote{
If $\sigma^2$ denotes the sample variance computed from $M = 10^6$ simulations,
the $95\%$ confidence interval for the variance is given by
$\left[ \sigma^2 \frac{M-1}{\chi^{-1}(1 - \frac{\alpha}{2},M-1)},
\sigma^2 \frac{M-1}{\chi^{-1}(\frac{\alpha}{2},M-1)} \right]$,
where $P = \chi(X,V)$ is the chi-square cumulative distribution function
with $V$ degrees of freedom \cite{Sh97}.
}.
In view of the result for the case that $\phi$ and $\psi$ are piecewise-constant
(Theorem \ref{th:yVarPWL}, given below),
we propose the following result:
\begin{conjecture}
Consider (\ref{eq:e}) with (\ref{eq:aLaRpos}).
Suppose $x(0)$ is random with \pdf, $p_{{\rm qss},\ee}$, and $y(0) = y_0$.
Then for any $\delta > 0$, whenever $\ee^{1-\delta} \le t \le \ee^{-M}$ we have
\begin{equation}
Var(y(t)) = \ee \kappa t + \frac{(b_L-b_R)^2}{(a_L+a_R)^2} \ee t + O(\ee^2) \;.
\label{eq:yVar2}
\end{equation}
\label{cj:yVar}
\end{conjecture}

From here until the concluding section, \S\ref{sec:CONC},
we study the case that $\phi$ and $\psi$ are piecewise-constant, i.e.
\begin{eqnarray}
\phi(x) &=& \left\{ \begin{array}{lc}
a_L \;, & x < 0 \\
-a_R \;, & x > 0
\end{array} \right. \;, \label{eq:tvphi} \\
\psi(x) &=& \left\{ \begin{array}{lc}
b_L \;, & x < 0 \\
b_R \;, & x > 0
\end{array} \right. \;. \label{eq:tvpsi}
\end{eqnarray}
\begin{theorem}
Consider (\ref{eq:e}) with (\ref{eq:aLaRpos})
and suppose $\phi$ and $\psi$ are given by (\ref{eq:tvphi})-(\ref{eq:tvpsi}).
Suppose $x(0)$ is random with \pdf, $p_{{\rm qss},\ee}$, and $y(0) = y_0$.
Then for any $\delta > 0$, whenever $t \ge \ee^{1-\delta}$ we have
\begin{equation}
Var(y(t)) = \ee \kappa t + \frac{(b_L-b_R)^2}{(a_L+a_R)^2} \ee t + O(\ee^2) \;.
\label{eq:yVarPWL}
\end{equation}
\label{th:yVarPWL}
\end{theorem}
We prove this result in the next section.

\section{Two-valued drift and a proof of Theorem \ref{th:yVarPWL}}
\label{sec:CALC}

The stochastic differential equation (\ref{eq:dx}) with (\ref{eq:tvphi}):
\begin{equation}
dx = \left\{ \begin{array}{lc}
a_L \;, & x<0 \\
-a_R \;, & x>0
\end{array} \right\} \, dt
+ \sqrt{\ee} \, dW_1(t) \;,
\label{eq:dxtv}
\end{equation}
has been referred to as Brownian motion with two-valued drift.
The transitional \pdf~of this process was first derived
by Karatzas and Shreve in \cite{KaSh84}.
In this section we state this \pdf~and use it prove Theorem \ref{th:yVarPWL}.

\subsection{The transitional probability density function for $x(t)$}

\begin{figure}[t!]
\begin{center}
\setlength{\unitlength}{1cm}
\begin{picture}(10,5)
\put(0,0){\includegraphics[height=5cm]{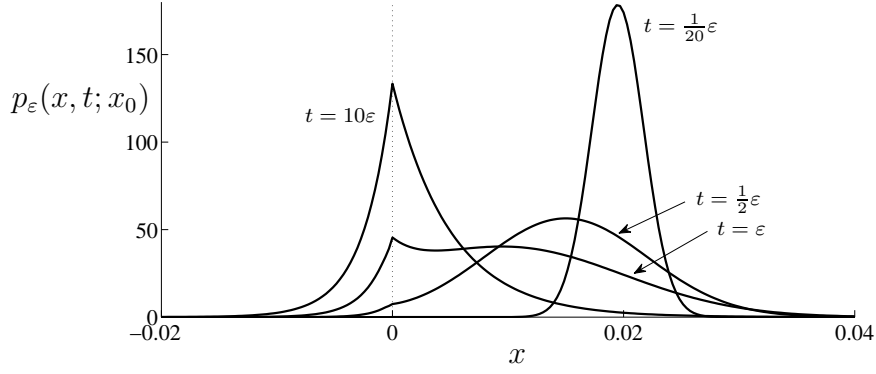}}
\put(5.1,.2){$x$}
\put(-1.5,3.6){$p_\ee(x,t;x_0)$}
\put(6.9,4.6){\scriptsize $t = \frac{1}{20} \ee$}
\put(7.6,2.3){\scriptsize $t = \frac{1}{2} \ee$}
\put(7.88,1.9){\scriptsize $t = \ee$}
\put(2.4,3.4){\scriptsize $t = 10 \ee$}
\end{picture}
\caption{
The probability density function of $x(t)$, (\ref{eq:p}), with
$x_0 = 0.02$,
$a_L = 2$,
$a_R = 1$ and
$\ee = 0.01$
at four different times.
\label{fig:examplePDFx}
}
\end{center}
\end{figure}

The transitional \pdf~for (\ref{eq:dxtv}) is given by
\begin{equation}
\hspace{-12mm}
p_\ee(x,t|x_0) = \left\{ \begin{array}{lc}
\frac{2}{\ee} {\rm e}^{\frac{2 a_L x}{\ee}} \int_0^\infty 
h_\ee(t,b,a_R) * h_\ee(t,b-x-x_0,a_L) \,db +
G_{{\rm absorb},\ee}(x,t,a_L|x_0) \;, & x_0 \le 0, \, x \le 0 \\
\frac{2}{\ee} {\rm e}^{\frac{-2 a_R x}{\ee}} \int_0^\infty
h_\ee(t,b+x,a_R) * h_\ee(t,b-x_0,a_L) \,db \;, & x_0 \le 0, \, x \ge 0 \\
\frac{2}{\ee} {\rm e}^{\frac{2 a_L x}{\ee}} \int_0^\infty
h_\ee(t,b+x_0,a_R) * h_\ee(t,b-x,a_L) \,db \;, & x_0 \ge 0, \, x \le 0 \\
\frac{2}{\ee} {\rm e}^{\frac{-2 a_R x}{\ee}} \int_0^\infty
h_\ee(t,b+x+x_0,a_R) * h_\ee(t,b,a_L) \,db +
G_{{\rm absorb},\ee}(x,t,-a_R|x_0) \;, & x_0 \ge 0, \, x \ge 0
\end{array} \right. \;,
\label{eq:p}
\end{equation}
where
\begin{equation}
h_\ee(t,x_0,\mu) \equiv \frac{|x_0|}{\sqrt{2 \pi \ee t^3}}
{\rm e}^{-\frac{(x_0 - \mu t)^2}{2 \ee t}} \;,
\label{eq:h}
\end{equation}
is the \pdf~for the first passage time to zero
of Brownian motion with constant drift\removableFootnote{
I've taken the definitions of $h$ and $G$ from \cite{KaSh91,Zh90},
and generalized to $\ee \ne 0$.
Both definitions involve $-\mu t$, which seems nice,
however the stochastic process associated with $h$ is
$dx = -\mu \, dt + \sqrt{\ee} \, dW(t)$,
whereas the stochastic process associated with $G$ and $G_{\rm absorb}$ is
$dx = \mu \, dt + \sqrt{\ee} \, dW(t)$.
},
\begin{equation}
G_{{\rm absorb},\ee}(x,t,\mu|x_0) \equiv
\frac{1}{\sqrt{2 \pi \ee t}}
{\rm e}^{-\frac{(x-x_0-\mu t)^2}{2 \ee t}} - {\rm e}^{\frac{-2 \mu x_0}{\ee}} 
\frac{1}{\sqrt{2 \pi \ee t}}
{\rm e}^{-\frac{(x+x_0-\mu t)^2}{2 \ee t}} \;.
\label{eq:Gabsorb}
\end{equation}
is the transitional \pdf~for
Brownian motion with constant drift
and an absorbing boundary condition at zero, and
\begin{equation}
f_1(t) * f_2(t) \equiv \int_0^t f_1(\tau) f_2(t-\tau) \, d\tau \;,
\end{equation}
is the convolution relating to Laplace transforms.
Fig.~\ref{fig:examplePDFx} shows (\ref{eq:p}) at different times\removableFootnote{
In the special case
\begin{equation}
a_L = a_R = a \;,
\end{equation}
we can use the convolution rule
\begin{equation}
h_\ee(t,z_1,\mu) * h_\ee(t,z_2,\mu) = \frac{1}{\ee} h_\ee(t,z_1+z_2,\mu) \;,
\label{eq:hconvrule}
\end{equation}
to simplify (\ref{eq:p}) to
\begin{equation}
p_\ee(x,t|x_0) = \frac{1}{\sqrt{2 \pi \ee t}}
{\rm e}^{\frac{a}{\ee} \left( |x_0| - |x| - \frac{at}{2} \right)}
{\rm e}^{\frac{-(x-x_0)^2}{2 \ee t}} +
\frac{a}{2 \ee}
{\rm e}^{\frac{-2 a |x|}{\ee}}
\left( 1 - {\rm erf} \left(
\frac{|x|+|x_0|-at}{\sqrt{2 \ee t}}
\right) \right) \;,
\label{eq:PDFxSym}
\end{equation}
where
\begin{equation}
{\rm erf}(z) = \frac{2}{\sqrt{\pi}} \int_0^z {\rm e}^{-u^2} \,du \;,
\label{eq:erf}
\end{equation}
is the error function.
}.

In \cite{KaSh84}, Karatzas and Shreve derive (\ref{eq:p})
for $\ee = 1$ by using Girsanov's theorem \cite{Gi60,Ok03,KaSh91}
and the trivariate \pdf~of Brownian motion,
its {\em positive occupation time}, and its {\em local time} about zero.
The result for $\ee \ne 1$ follows simply from the scaling (\ref{eq:xthat}).
The Laplace transform of (\ref{eq:p}) can be written as an integral-free
expression and this was achieved in the earlier paper \cite{BeSh80}.
We do this below for $x_0 = 0$ and $x > 0$.
In \cite{QiZh02,QiRu03}, the authors derive $p_\ee(0,t|x_0)$
and use this to bound \pdf s for a large class
of scalar stochastic differential equations.
In \cite{GrHe01}, $p_\ee(x,t|0)$ is studied in the case $a_L, a_R < 0$.
In \cite{Zh90}, Zhang derived an expression for the transitional \pdf~of
Brownian motion with a general bounded piecewise-continuous drift function.
This could be used to analyze the \pdf~of (\ref{eq:dx}) with general $\phi$ asymptotically,
but such a calculation is beyond of scope of this paper.

\subsection{Proof of Theorem \ref{th:yVarPWL}}

In the case of two-valued drift (\ref{eq:tvphi}),
the quasi-steady-state density (\ref{eq:qss}) is a true steady-state
defined for all $x \in \mathbb{R}$:
\begin{equation}
p_{{\rm ss},\ee}(x) = \left\{ \begin{array}{lc}
\frac{K}{\ee} {\rm e}^{\frac{2 a_L x}{\ee}} \;, & x<0 \\
\frac{K}{\ee} {\rm e}^{\frac{-2 a_R x}{\ee}} \;, & x>0
\end{array} \right. \;, \qquad
K = \frac{2 a_L a_R}{a_L+a_R} \;.
\label{eq:ss}
\end{equation}
With the notation
\begin{equation}
\frac{\partial p_\ee^\pm}{\partial x}(0,t|x_0)
\equiv \lim_{\Delta \to 0^\pm} \frac{\partial p_\ee}{\partial x}(\Delta,t|x_0) \;,
\end{equation}
we have the following expression for the probability that $x(t) > 0$,
given $x_0 = 0$.

\begin{lemma}
For any $t > 0$,
\begin{equation}
\int_0^\infty p_\ee(x,t|0) \,dx =
\frac{1}{2} -
\int_0^t a_R p_\ee(0,s|0) +
\frac{\ee}{2} \frac{\partial p_\ee^+}{\partial x}(0,s|0) \, ds \;.
\label{eq:intPDFpos}
\end{equation}
\label{le:intPDFpos}
\end{lemma}

\myproof
For $x>0$, in the case of two-valued drift (\ref{eq:tvphi}),
the Fokker-Planck equation (\ref{eq:fp}) is
\begin{equation}
\frac{\partial p_\ee}{\partial t} = \frac{\partial}{\partial x} \left(
a_R p_\ee + \frac{\ee}{2} \frac{\partial p_\ee}{\partial x} \right) \;.
\end{equation}
Integration over $x \in [\Delta,\infty)$ yields
\begin{equation}
\frac{\partial}{\partial t} \int_\Delta^\infty p_\ee(x,t|0) \,dx =
-\left( a_R p_\ee(\Delta,t|0) + \frac{\ee}{2} \frac{\partial p_\ee}{\partial x}(\Delta,s|0) \right) \;,
\label{eq:FPintegrated}
\end{equation}
for any $\Delta > 0$.
Then integrating with respect to $t$ and taking $\Delta \to 0$
produces (\ref{eq:intPDFpos}) where we also use\removableFootnote{
The identity (\ref{eq:intPDF0}) is proved for (\ref{eq:dx})
for any bounded $\ee$-independent function $\phi$, by noting that as $\ee \to 0^+$,
$p_\ee(x,t|0)$ approaches a zero-mean Gaussian.
Indeed the scaling
\begin{equation}
\tilde{x} = \frac{x}{\ee^2} \;, \qquad
\tilde{t} = \frac{t}{\ee^3} \;,
\end{equation}
moves the $\ee$-dependence in (\ref{eq:dx}) from the diffusion term to the drift term:
\begin{equation}
d\tilde{x} = \ee \phi(\ee^2 \tilde{x}) \,d\tilde{t} + dW_1(\tilde{t}) \;.
\end{equation}
Then the \pdf~for this stochastic differential equation, 
admits a regular expansion in $\ee$ and thus
\begin{equation}
p(x,t|0) = \frac{1}{\ee^2}
\left( \frac{1}{2 \pi \tilde{t}} \,{\rm e}^{-\frac{\tilde{x}^2}{2 \tilde{t}}} + O(\ee) \right) \;,
\end{equation}
and the result follows.
}
\begin{equation}
\lim_{t \to 0^+} \int_0^\infty p_\ee(x,t|0) \,dx = \frac{1}{2} \;,
\label{eq:intPDF0}
\end{equation}
which may be demonstrated by noting that as $t \to 0^+$,
$p_\ee(x,t|0)$ is well-approximated by a zero-mean Gaussian.
\hfill
$\Box$

Proofs of the following two lemmas are given in Appendix \ref{sec:PROOFS}.
Theorem \ref{th:yVarPWL} is an immediate consequence of
Lemma \ref{le:sgnxssgnxu} combined with
(\ref{eq:sgnxMean}) and (\ref{eq:varIdentity}).

\begin{lemma}
For the density (\ref{eq:p}),
\begin{eqnarray}
\int_0^\infty a_R p_\ee(0,t|0) + \frac{\ee}{2} \frac{\partial p_\ee^+}{\partial x}(0,t|0) \,dt &=&
\frac{-(a_L-a_R)}{2 (a_L+a_R)} \;, \label{eq:int1} \\
\int_0^\infty t \left( a_R p_\ee(0,t|0) + \frac{\ee}{2} \frac{\partial p_\ee^+}{\partial x}(0,t|0) \right) \,dt &=&
\frac{-\ee (a_L-a_R)}{2 a_L a_R (a_L+a_R)} \;. \label{eq:int2}
\end{eqnarray}
\label{le:ints}
\end{lemma}

\begin{lemma}
Consider (\ref{eq:dxtv}) and suppose $a_L, a_R > 0$ and
$x(0)$ is random with \pdf, $p_{{\rm ss},\ee}$.
Then for any $\delta > 0$, if $t \ge \ee^{1-\delta}$, we have
\begin{equation}
\int_0^t \int_0^t \big< {\rm sgn}(x(s)) \, {\rm sgn}(x(u)) \big> \, ds \, du =
\frac{(a_L-a_R)^2 t^2}{(a_L+a_R)^2} + \frac{4 \ee t}{(a_L+a_R)^2} + O(\ee^2) \;.
\label{eq:sgnxssgnxu}
\end{equation}
\label{le:sgnxssgnxu}
\end{lemma}

In the special case, $a_L = a_R = a$, we can write $Var(y(t))$ exactly.
In this case, by symmetry, $\int_0^\infty p_\ee(x,t|0) \,dx \equiv \frac{1}{2}$,
thus by Lemma \ref{le:intPDFpos}, 
$a p_\ee(0,t|0) + \frac{\ee}{2} \frac{\partial p_\ee^+}{\partial x}(0,t|0) \equiv 0$.
Consequently the integral on the right-hand side of (\ref{eq:sgnxssgnxu4}) vanishes,
and (\ref{eq:Q2}) leads to
\begin{eqnarray}
{\rm Var}(y(t)) &=& \frac{(b_L-b_R)^2}{4} \Bigg(
\frac{\ee t}{a^2} - \frac{\ee^2}{a^4} +
\frac{\sqrt{2t}}{\sqrt{\pi \ee}} \left(
\frac{\ee^2}{a^3} - \frac{2 \ee t}{3 a} - \frac{a t^2}{3} \right)
{\rm e}^{-\frac{a^2 t}{2 \ee}} \nonumber \\
&&+~\left( \frac{\ee^2}{a^4} - \frac{\ee t}{a^2} + t^2 + \frac{a^2 t^3}{3 \ee} \right)
{\rm erfc} \left( \frac{a \sqrt{t}}{\sqrt{2 \ee}} \right) \Bigg) \;,
\label{eq:yVarPWLSym}
\end{eqnarray}
which is consistent with (\ref{eq:yVarPWL}).
\section{Conclusions}
\label{sec:CONC}

When small noise is added to a Filippov system,
orbits no longer slide along an attracting sliding section of a switching manifold.
Instead, with high probability, orbits follow a random path near the switching manifold,
Fig.~\ref{fig:noisyFlowExSlide}-B.
The average size of deviations from the switching manifold is governed by the strength
of the noise relative to the magnitude of the vector field
in a direction orthogonal to the switching manifold.

The general $N$-dimensional stochastic differential equation
formed by adding noise to the Filippov system (\ref{eq:Filippov2})
is particularly difficult to analyze due to the discontinuity and multiple dimensions.
For this reason we made the supposition
that the system is invariant along the switching manifold.
This prevents an exploration of the effects of noise on orbits
that reach the end of an attracting sliding region,
but enables calculations to be reduced to one dimension.
Moreover, with this reduction directions parallel to the switching manifold
may be treated identically and hence it is sufficient
to study the two-dimensional system (\ref{eq:e}).

For the system (\ref{eq:e}), $x(t)$ denotes the displacement from the switching manifold
and is governed by (\ref{eq:dx}) with (\ref{eq:phi}).
Sample paths of (\ref{eq:dx}) settle to the quasi-steady-state \pdf,
$p_{{\rm qss},\ee}$ (\ref{eq:qss}),
on an $O(\ee)$ time scale.
This \pdf~is not a true steady-state because orbits escape a neighbourhood of zero.
However, since this occurs on an exponentially long time scale,
see \S\ref{sub:EXIT}, it is suitable to assume $x(t)$
is distributed by $p_{{\rm qss},\ee}$ at times in a given range
$\ee^{1-\delta} \le t \le \ee^{-M}$, for any $\delta, M > 0$,
where we may take small $\delta$ and large $M$, such that this is a long time interval.
This assumption has the benefit of significantly simplifying our calculations.

Equation (\ref{eq:dy}) is the stochastic differential equation
for $y(t)$, which represents displacement along the switching manifold.
In the limit, $\ee \to 0$,
the mean of $y(t)$, (\ref{eq:yMean}), limits on Filippov's sliding solution, $y_{\rm slide}(t)$.
For $\ee \ne 0$, the perturbation of $\big< y(t) \big>$ from $y_{\rm slide}(t)$ is $O(\ee)$,
as for a generic smooth system.
The perturbation depends on linear terms in $\phi$ and $\psi$.
In order to gauge the effect of this perturbation on overall dynamics,
it is necessary to compare it to the standard deviation of $y(t)$.
In the case that $\phi$ and $\psi$ are piecewise-constant,
${\rm Var}(y(t))$ is $O(\ee)$, see Theorem \ref{th:yVarPWL}.
For general $\phi$ and $\psi$,
we conjectured that the leading order terms of ${\rm Var}(y(t))$ are unchanged.
Consequently deviations of $y(t)$ from $\big< y(t) \big>$ are $O(\sqrt{\ee})$.
Therefore, assuming $\ee \ll 1$, we expect these deviations to dominate the
difference between $\big< y(t) \big>$ and $y_{\rm slide}(t)$.

Although our above calculations are for differential equations
that are static along the switching manifold,
we can apply the basic principles gained to more general systems
such as the relay control system, (\ref{eq:relayControlSystem}) with (\ref{eq:ABCvalues}).
First note that since the deterministic equations are independent of $\ee$,
when we make statements asymptotically in $\ee$,
we implicitly assume that the magnitude of 
any parameter in the deterministic equations is much less than,
in particular, $\frac{1}{\sqrt{\ee}}$.
But, for (\ref{eq:relayControlSystem}) with (\ref{eq:ABCvalues}),
$\frac{\partial \dot{x_2}}{\partial x_1} = -\zeta-100$,
which is an extremely large value (we used $\zeta = 0.06$).
Since $x_1$ represents displacement from the switching manifold
and $x_2$ is a direction parallel to the switching manifold,
in the context of (\ref{eq:e})-(\ref{eq:psi}),
$\frac{\partial \dot{x_2}}{\partial x_1}$ corresponds to the values of $d_L$ and $d_R$.
These values influence, to lowest order, the perturbation of the mean from the
deterministic solution, (\ref{eq:yMean}),
but not the deviation of sample paths from the mean, (\ref{eq:yVar2}).
This is a possible explanation for the noise-induced effect
identified in \S\ref{sec:MOTIV}\removableFootnote{
This would predict the mean oscillation time to differ from the deterministic period by $O(\ee)$,
and the standard deviation of the oscillation time to be $O(\sqrt{\ee})$.
This seems to disagree with Fig.~\ref{fig:noisyExample}-B,
but I think this is because the figure involves relatively large values of $\ee$.
The scaling predictions may be justified (somewhat)
by a plot involving $\ee$ on the horizontal axis and $\ee_{\rm max} \approx 0.0003$.
}.
Specifically we found that when the parameters of
(\ref{eq:relayControlSystem}) with (\ref{eq:ABCvalues})
are tuned such that in the absence of noise orbits settle
to oscillatory motion with sliding,
the addition of noise may significantly decrease
the average oscillation time, Fig.~\ref{fig:noisyExample}.
A further analysis is required in order to make more definitive statements.
We believe that the decrease in oscillation time with $\ee$ is
a consequence of particular geometrical orientations and that
an increase in oscillation time with $\ee$ is equally possible for this type of system.

This paper leaves many avenues for future investigations.
Perhaps foremost, we would like to understand the perturbation of $\big< y(t) \big>$
from $y_{\rm slide}$ in the case that the system explicitly depends on $y$.
In particular we would like to understand how noise effects dynamics
near the end of an attracting sliding region.
Also it remains to consider more general forms for the noise,
such as coloured noise or noise that is correlated in $x$ and $y$,
and study the effects of noise on other scenarios such as sliding bifurcations \cite{CoDi12}.

\appendix

\section{Proofs}
\label{sec:PROOFS}

\subsection{Proof of Lemma \ref{le:ints}}
\label{sub:PROOF1}

By Lemma \ref{le:intPDFpos} and (\ref{eq:ss}),
\begin{equation}
\int_0^\infty a_R p_\ee(0,t|0) + \frac{\ee}{2} \frac{\partial p_\ee^+}{\partial x}(0,t|0) \,dt
= \frac{1}{2} - \int_0^\infty p_{{\rm ss},\ee}(x) \,dx
= \frac{-(a_L-a_R)}{2 (a_L+a_R)} \;,
\end{equation}
which verifies (\ref{eq:int1}).
To obtain (\ref{eq:int2}), we multiply (\ref{eq:FPintegrated}) by $t$ and
take $\Delta \to 0$ to obtain
\begin{equation}
t \frac{\partial}{\partial t} \int_0^\infty p_\ee(x,t|0) \,dx =
-t \left( a_R p_\ee(0,t|0) + \frac{\ee}{2} \frac{\partial p_\ee^+}{\partial x}(0,t|0) \right) \;.
\end{equation}
Integration by parts yields
\begin{equation}
\int_0^T t \left( a_R p_\ee(0,t|0) + \frac{\ee}{2} \frac{\partial p_\ee^+}{\partial x}(0,t|0) \right) \,dt
= \int_0^T \int_0^\infty p_\ee(x,t|0) - p_\ee(x,T|0) \,dx \,dt \;,
\end{equation}
for any $T>0$.
Taking $T \to \infty$ gives
\begin{equation}
\int_0^\infty t \left( a_R p_\ee(0,t|0) + \frac{\ee}{2} \frac{\partial p_\ee^+}{\partial x}(0,t|0) \right) \,dt =
\int_0^\infty \int_0^\infty p_\ee(x,t|0) - p_{{\rm ss},\ee}(x) \,dx \,dt \;.
\label{eq:int2b}
\end{equation}
We evaluate the double integral using Laplace transforms.
For $x>0$,
\begin{equation}
P_{{\rm ss},\ee}(x,\lambda) \equiv \mathcal{L}[ p_{{\rm ss},\ee}(x) ] =
\int_0^\infty {\rm e}^{-\lambda t} K {\rm e}^{\frac{-2 a_R x}{\ee}} \,dt =
\frac{2 a_L a_R {\rm e}^{\frac{-2 a_R x}{\ee}}}{\ee (a_L + a_R) \lambda} \;.
\label{eq:Pss}
\end{equation}
The Laplace transform of $h_\ee$, (\ref{eq:h}), is
\begin{equation}
H_\ee(\lambda,z,\mu) \equiv \mathcal{L}[ h_\ee(t,z,\mu) ] =
{\rm e}^{\frac{1}{\ee} \left( \mu z - \sqrt{\mu^2 + 2 \ee \lambda} |z| \right)} \;,
\label{eq:H}
\end{equation}
for $\lambda > 0$.
Using (\ref{eq:p}) and expanding about $\lambda = 0$:
\begin{eqnarray}
P_\ee(x,\lambda|0)
&\equiv& \mathcal{L}[ p_\ee(x,t|0) ] \nonumber \\
&=& \frac{2 {\rm e}^{-\frac{2 a_R x}{\ee}}}{\ee}
\int_0^\infty
H_\ee(\lambda,b+x,a_R) H_\ee(\lambda,b,a_L) \,db \nonumber \\
&=& 
\frac{2 {\rm e}^{-\frac{1}{\ee} \left( a_R +
\sqrt{a_R^2 + 2 \ee \lambda} \right) x}}
{-a_R + \sqrt{a_R^2 + 2 \ee \lambda} - a_L + \sqrt{a_L^2 + 2 \ee \lambda}} \nonumber \\
&=& \frac{2 {\rm e}^{-\frac{2 a_R x}{\ee}}}{\ee}
\left( \frac{a_L a_R}{(a_L+a_R) \lambda} +
\frac{\ee (a_L^3+a_R^3) - 2 a_L^2 a_R (a_L+a_R) x}
{2 a_L a_R (a_L+a_R)^2} + O(\lambda) \right) \;.
\label{eq:P}
\end{eqnarray}
By (\ref{eq:Pss}) and (\ref{eq:P}),
\begin{equation}
P_\ee(x,\lambda|0) - P_{{\rm ss},\ee}(x,\lambda) =
\frac{\ee (a_L^3+a_R^3) - 2 a_L^2 a_R (a_L+a_R) x}
{\ee a_L a_R (a_L+a_R)^2}
{\rm e}^{\frac{-2 a_R x}{\ee}} + O(\lambda) \;,
\end{equation}
and thus from (\ref{eq:int2b})
\begin{eqnarray}
\int_0^\infty t \left( a_R p_\ee(0,t|0) + \frac{\ee}{2} \frac{\partial p_\ee^+}{\partial x}(0,t|0) \right) \,dt
&=& \int_0^\infty P_\ee(x,0|0) - P_{{\rm ss},\ee}(x,0) \,dx \nonumber \\
&=& \int_0^\infty 
\frac{\ee (a_L^3+a_R^3) - 2 a_L^2 a_R (a_L+a_R) x}
{\ee a_L a_R (a_L+a_R)^2}
{\rm e}^{\frac{-2 a_R x}{\ee}} \,dx \nonumber \\
&=& \frac{-\ee (a_L-a_R)}{2 a_L a_R (a_L+a_R)} \;,
\end{eqnarray}
as required. \hfill $\Box$

\subsection{Proof of Lemma \ref{le:sgnxssgnxu}}
\label{sub:PROOF2}

Since (\ref{eq:dx}) has no explicit time-dependence,
\begin{eqnarray}
\int_0^t \int_0^t \big< {\rm sgn}(x(s)) \, {\rm sgn}(x(u)) \big> \,ds \,du
&=& \int_0^t \int_0^t \big< {\rm sgn}(x(0)) \, {\rm sgn}(x(|u-s|)) \big> \,ds \,du \nonumber \\
&=& 2 \int_0^t (t-u) \big< {\rm sgn}(x(0)) \, {\rm sgn}(x(u)) \big> \,du \;.
\label{eq:sgnxssgnxu2}
\end{eqnarray}
Furthermore,
\begin{eqnarray}
\big< {\rm sgn}(x(0)) \, {\rm sgn}(x(u)) \big>
&=& \int_{-\infty}^\infty \int_{-\infty}^\infty
{\rm sgn}(x) p_{{\rm ss},\ee}(x) \, {\rm sgn}(y) p_\ee(y,u|x) \,dy \,dx \nonumber \\
&=& -\int_{-\infty}^0 K {\rm e}^{\frac{2 a_L x}{\ee}}
\int_{-\infty}^\infty {\rm sgn}(y) p_\ee(y,u|x) \,dy \,dx \nonumber \\
&&+~\int_0^\infty K {\rm e}^{\frac{-2 a_R x}{\ee}}
\int_{-\infty}^\infty {\rm sgn}(y) p_\ee(y,u|x) \,dy \,dx \;.
\label{eq:sgnx0sgnxu}
\end{eqnarray}
To evaluate (\ref{eq:sgnxssgnxu2}) using (\ref{eq:sgnx0sgnxu}),
we reorder the integrals of $y$, $x$ and $u$,
but to do this we must first transfer the $u$-dependence
from the integrand to the limits of integration.
We have
\begin{eqnarray}
\int_{-\infty}^\infty {\rm sgn}(y) p_\ee(y,u|x<0) \,dy
&=& 2 \int_0^\infty p_\ee(y,u|x<0) \,dy - 1 \nonumber \\
&=& 2 \left( \int_0^u h_\ee(s,x,-a_L)
\int_0^\infty p_\ee(y,u-s|0) \,dy \,ds \right) - 1 \;,
\end{eqnarray}
where in the second line
we conditioned over the first passage time, $s$, of (\ref{eq:dx})
from $x<0$ to $0$.
By Lemma \ref{le:intPDFpos}, and since
$\int_0^\infty h_\ee(s,x,-a_L) \,ds = 1$,
for all $x<0$, we obtain
\begin{eqnarray}
\int_{-\infty}^\infty {\rm sgn}(y) p_\ee(y,u|x<0) \,dy
&=&
-\int_u^\infty h_\ee(s,x,-a_L) \,ds \nonumber \\
&&-~2 \int_0^u h_\ee(s,x,-a_L)
\int_0^{u-s} a_R p_\ee(0,v|0) + \frac{\ee}{2} \frac{\partial p_\ee^+}{\partial x}(0,v|0) \,dv \,ds \;.
\label{eq:sgnypxneg}
\end{eqnarray}
A similar calculation with $x>0$ produces
\begin{eqnarray}
\int_{-\infty}^\infty {\rm sgn}(y) p_\ee(y,u|x>0) \,dy
&=&
\int_u^\infty h_\ee(s,x,a_R) \,ds \nonumber \\
&&-~2 \int_0^u h_\ee(s,x,a_R)
\int_0^{u-s} a_R p_\ee(0,v|0) + \frac{\ee}{2} \frac{\partial p_\ee^+}{\partial x}(0,v|0) \,dv \,ds \;.
\label{eq:sgnypxpos}
\end{eqnarray}
By combining (\ref{eq:sgnxssgnxu2}), (\ref{eq:sgnx0sgnxu}),
(\ref{eq:sgnypxneg}) and (\ref{eq:sgnypxpos}) we arrive at
\begin{eqnarray}
&& \int_0^t \int_0^t \big< {\rm sgn}(x(s)) \, {\rm sgn}(x(u)) \big> \,ds \,du
= 2 \int_0^t (t-u) \int_{-\infty}^0 K {\rm e}^{\frac{2 a_L x}{\ee}}
\int_u^\infty h_\ee(s,x,-a_L)
\,ds \,dx \,du \nonumber \\
&&\hspace{10mm}+~4 \int_0^t (t-u) \int_{-\infty}^0 K {\rm e}^{\frac{2 a_L x}{\ee}}
\int_0^u h_\ee(s,x,-a_L)
\int_0^{u-s} a_R p_\ee(0,v|0) + \frac{\ee}{2} \frac{\partial p_\ee^+}{\partial x}(0,v|0)
\,dv \,ds \,dx \,du \nonumber \\
&&+~2 \int_0^t (t-u) \int_0^\infty K {\rm e}^{\frac{-2 a_R x}{\ee}}
\int_u^\infty h_\ee(s,x,a_R)
\,ds \,dx \,du \nonumber \\
&&\hspace{-10mm}-~4 \int_0^t (t-u) \int_0^\infty K {\rm e}^{\frac{-2 a_R x}{\ee}}
\int_0^u h_\ee(s,x,a_R)
\int_0^{u-s} a_R p_\ee(0,v|0) + \frac{\ee}{2} \frac{\partial p_\ee^+}{\partial x}(0,v|0)
\,dv \,ds \,dx \,du \;.
\label{eq:sgnxssgnxu3}
\end{eqnarray}
We now simplify the four terms of (\ref{eq:sgnxssgnxu3}).
We define
\begin{equation}
Q_\ee(t,a) \equiv \int_0^t (t-u) \int_0^\infty K {\rm e}^{\frac{-2 a x}{\ee}}
\int_0^u h_\ee(s,x,a) \,ds \,dx \,du \;,
\label{eq:Q}
\end{equation}
so that the first term of (\ref{eq:sgnxssgnxu3}) may be written as
\begin{equation}
2 \int_0^t (t-u) \int_{-\infty}^0 K {\rm e}^{\frac{2 a_L x}{\ee}}
\int_u^\infty h_\ee(s,x,-a_L)
\,ds \,dx \,du =
\frac{a_R}{a_L+a_R} t^2 - 2 Q_\ee(t,a_L) \;.
\end{equation}
Similarly the third term is
\begin{equation}
2 \int_0^t (t-u) \int_0^\infty K {\rm e}^{\frac{-2 a_R x}{\ee}} \,dx \,du -
2 Q_\ee(t,a_R) =
\frac{a_L}{a_L+a_R} t^2 - 2 Q_\ee(t,a_R) \;.
\end{equation}
To the second term of (\ref{eq:sgnxssgnxu3})
we reorder the integration such that $dv$ is the outer-most integral
instead of the inner-most integral.
This step is straight-forward but requires some care with
the limits of integration.
We obtain
\begin{eqnarray}
4 \int_0^t (t-u) \int_{-\infty}^0 K {\rm e}^{\frac{2 a_L x}{\ee}}
\int_0^u h_\ee(s,x,-a_L)
\int_0^{u-s} a_R p_\ee(0,v|0) + \frac{\ee}{2}
\frac{\partial p_\ee^+}{\partial x}(0,v|0)
\,dv \,ds \,dx \,du \nonumber \\
= 4 \int_0^t \left( a_R p_\ee(0,v|0) + \frac{\ee}{2}
\frac{\partial p_\ee^+}{\partial x}(0,v|0) \right) Q_\ee(t-v,a_L) \,dv \;,
\end{eqnarray}
and similarly for the last term of (\ref{eq:sgnxssgnxu3}):
\begin{eqnarray}
-4 \int_0^t (t-u) \int_0^\infty K {\rm e}^{\frac{-2 a_R x}{\ee}}
\int_0^u h_\ee(s,x,a_R)
\int_0^{u-s} a_R p_\ee(0,v|0) + \frac{\ee}{2} \frac{\partial p_\ee^+}{\partial x}(0,v|0)
\,dv \,ds \,dx \,du \nonumber \\
= -4 \int_0^t \left( a_R p_\ee(0,v|0) + \frac{\ee}{2} \frac{\partial p_\ee^+}{\partial x}(0,v|0) \right) Q_\ee(t-v,a_R) \,dv \;.
\end{eqnarray}
We are now able to write (\ref{eq:sgnxssgnxu3}) as
\begin{eqnarray}
&& \int_0^t \int_0^t \big< {\rm sgn}(x(s)) \, {\rm sgn}(x(u)) \big> \,ds \,du = 
t^2 - 2 \big( Q_\ee(t,a_L) + Q_\ee(t,a_R) \big) \nonumber \\
&&\hspace{10mm}+~4 \int_0^t \left( a_R p_\ee(0,v|0) + \frac{\ee}{2} \frac{\partial p_\ee^+}{\partial x}(0,v|0) \right)
\big( Q_\ee(t-v,a_L) - Q_\ee(t-v,a_R) \big) \,dv \;.
\label{eq:sgnxssgnxu4}
\end{eqnarray}
Via (\ref{eq:p}),
it may be demonstrated that
$a_R p_\ee(0,v|0) + \frac{\ee}{2} \frac{\partial p_\ee^+}{\partial x}(0,v|0)$
decays exponentially to zero as $v \to \infty$
on an $O(\ee)$ time scale.
For this reason it is helpful to 
apply the substitution $\tilde{v} = \frac{v}{\ee}$ to obtain
\begin{eqnarray}
&& \int_0^t \int_0^t \big< {\rm sgn}(x(s)) \, {\rm sgn}(x(u)) \big> \,ds \,du = 
t^2 - 2 \big( Q_\ee(t,a_L) + Q_\ee(t,a_R) \big) \nonumber \\
&&\hspace{10mm}+~4 \ee \int_0^{\frac{t}{\ee}} \left(
a_R p_\ee(0,\ee \tilde{v}|0) + \frac{\ee}{2} \frac{\partial p_\ee^+}{\partial x}(0,\ee \tilde{v}|0) \right)
\big( Q_\ee(t-\ee \tilde{v},a_L) - Q_\ee(t-\ee \tilde{v},a_R) \big) \,d\tilde{v} \;,
\label{eq:sgnxssgnxu5}
\end{eqnarray}
and expand
$Q_\ee(t-\ee \tilde{v},a_L) - Q_\ee(t-\ee \tilde{v},a_R)$ in $\ee$
such that the integral on the right-hand side of (\ref{eq:sgnxssgnxu5}) has the form
\begin{equation}
\sum_i \sum_j \alpha_{ij} \ee^i
\int_0^{\frac{t}{\ee}} \tilde{v}^j
\left( a_R p_\ee(0,\ee \tilde{v}|0) + \frac{\ee}{2} \frac{\partial p_\ee^+}{\partial x}(0,\ee \tilde{v}|0) \right)
\,d\tilde{v} \;,
\label{eq:intExpansion}
\end{equation}
for some coefficients $\alpha_{ij}$.
To obtain the coefficients,
we first evaluate (\ref{eq:Q}) via
multiple applications of integration by parts:
\begin{eqnarray}
Q_\ee(t,a)
&=& K \left( \frac{\ee t^2}{4a} - \frac{\ee^2 t}{4a^3} + \frac{\ee^3}{4a^5} \right) +
\frac{K \sqrt{t}}{\sqrt{2\pi}} \left(
\frac{\sqrt{\ee} t^2}{6} + \frac{\ee^\frac{3}{2} t}{3a^2} - \frac{\ee^\frac{5}{2}}{2a^4} \right)
{\rm e}^{\frac{-a^2 t}{2 \ee}} \nonumber \\
&&-~K \left( \frac{a t^3}{12} + \frac{\ee t^2}{4a} -
\frac{\ee^2 t}{4a^3} + \frac{\ee^3}{4a^5} \right)
{\rm erfc} \left( \frac{a \sqrt{t}}{\sqrt{2 \ee}} \right)
\label{eq:Q2} \;.
\end{eqnarray}
In view of (\ref{eq:sgnxssgnxu5}), we use (\ref{eq:Q2}) to obtain
\begin{eqnarray}
Q_\ee(t,a_L) + Q_\ee(t,a_R)
&=& \frac{t^2}{2} -
\frac{(a_L^3+a_R^3) \ee t}{2 a_L^2 a_R^2 (a_L+a_R)} + O(\ee^2) \label{eq:QplusQ} \;, \\
Q_\ee(t - \ee \tilde{v},a_L) - Q_\ee(t - \ee \tilde{v},a_R)
&=& \frac{-(a_L-a_R) t^2}{2(a_L+a_R)} +
\frac{(a_L^3-a_R^3) \ee t}{2 a_L^2 a_R^2 (a_L+a_R)} +
\frac{(a_L-a_R) \ee t \tilde{v}}{(a_L+a_R)} + O(\ee^2) \label{eq:QminusQ} \;.
\end{eqnarray}
Using Laplace's method \cite{BeOr78} to asymptotically evaluate (\ref{eq:intExpansion}),
the substitution of (\ref{eq:QplusQ}) and (\ref{eq:QminusQ}) into (\ref{eq:sgnxssgnxu5})
yields
\begin{eqnarray}
&& \int_0^t \int_0^t \big< {\rm sgn}(x(s)) \, {\rm sgn}(x(u)) \big> \,ds \,du =
\frac{(a_L^3+a_R^3) \ee t}{a_L^2 a_R^2 (a_L+a_R)} \nonumber \\
&&\hspace{10mm}-~\left( \frac{2 (a_L-a_R) t^2}{(a_L+a_R)} +
\frac{2 (a_L^3-a_R^3) \ee t}{a_L^2 a_R^2 (a_L+a_R)} \right)
\int_0^\infty \left( a_R p_\ee(0,v|0) + \frac{\ee}{2} \frac{\partial p_\ee^+}{\partial x}(0,v|0) \right) \,dv \nonumber \\
&&\hspace{10mm}+~\frac{4 (a_L-a_R) t}{(a_L+a_R)}
\int_0^\infty v \left( a_R p_\ee(0,v|0) + \frac{\ee}{2} \frac{\partial p_\ee^+}{\partial x}(0,v|0) \right) \,dv
+ O(\ee^2) \;.
\end{eqnarray}
By applying Lemma \ref{le:ints} and simplifying we finally
arrive at (\ref{eq:sgnxssgnxu}). \hfill
$\Box$


\begin{thebibliography}{10}

\bibitem{DiBu08}
M.~di~Bernardo, C.J. Budd, A.R. Champneys, and P.~Kowalczyk.
\newblock {\em Piecewise-smooth Dynamical Systems. Theory and Applications.}
\newblock Springer-Verlag, New York, 2008.

\bibitem{VaSc00}
A.J. Van~der Schaft and J.M. Schumacher.
\newblock {\em An Introduction to Hybrid Dynamical Systems.}
\newblock Springer-Verlag, New York, 2000.

\bibitem{LeNi04}
R.I. Leine and H.~Nijmeijer.
\newblock {\em Dynamics and Bifurcations of Non-smooth Mechanical Systems},
  volume~18 of {\em Lecture Notes in Applied and Computational Mathematics}.
\newblock Springer-Verlag, Berlin, 2004.

\bibitem{BaVe01}
S.~Banerjee and G.C. Verghese, editors.
\newblock {\em Nonlinear Phenomena in Power Electronics.}
\newblock IEEE Press, New York, 2001.

\bibitem{ZhMo03}
Z.T. Zhusubaliyev and E.~Mosekilde.
\newblock {\em Bifurcations and Chaos in Piecewise-Smooth Dynamical Systems.}
\newblock World Scientific, Singapore, 2003.

\bibitem{PuSu06}
T.~Puu and I.~Sushko, editors.
\newblock {\em Business Cycle Dynamics: Models and Tools.}
\newblock Springer-Verlag, New York, 2006.

\bibitem{DaNo00}
H.~Dankowicz and A.B. Nordmark.
\newblock On the origin and bifurcations of stick-slip oscillations.
\newblock {\em Phys. D}, 136:280--302, 2000.

\bibitem{FrNo00}
M.H. Fredriksson and A.B. Nordmark.
\newblock On normal form calculation in impact oscillators.
\newblock {\em Proc. R. Soc. A}, 456:315--329, 2000.

\bibitem{DiBu01}
M.~di~Bernardo, C.J. Budd, and A.R. Champneys.
\newblock Normal form maps for grazing bifurcations in $n$-dimensional
  piecewise-smooth dynamical systems.
\newblock {\em Phys. D}, 160:222--254, 2001.

\bibitem{BeGe06}
N.~Berglund and B.~Gentz.
\newblock {\em Noise-Induced Phenomena in Slow-Fast Dynamical Systems.}
\newblock Springer, New York, 2006.

\bibitem{LiGa04}
B.~Lindner, J.~Garcia-Ojalvo, A.~Neiman, and L.~Schimansky-Geier.
\newblock Effects of noise in excitable systems.
\newblock {\em Phys. Reports}, 392:321--424, 2004.

\bibitem{PiKu97}
A.S. Pikovsky and J.~Kurths.
\newblock Coherence resonance in a noise-driven excitable system.
\newblock {\em Phys. Rev. Lett.}, 78(5):775--778, 1997.

\bibitem{GaHa98}
L.~Gammaitoni, P.~H\"{a}nggi, P.~Jung, and F.~Marchesoni.
\newblock Stochastic resonance.
\newblock {\em Rev. Modern Phys.}, 70(1):223--287, 1998.

\bibitem{Gr05}
T.C.L. Griffin.
\newblock {\em Dynamics of Stochastic Nonsmooth Systems.}
\newblock PhD thesis, University of Bristol, 2005.

\bibitem{Wa98}
R.~Wackerbauer.
\newblock Noise-induced stabilization of one-dimensional discontinuous maps.
\newblock {\em Phys. Rev. E}, 58(3):3036--3044, 1998.

\bibitem{ZhSh06}
L.~Zhang, P.~Shi, C.~Wang, and H.~Gao.
\newblock Robust {$H_\infty$} filtering for switched linear discrete-time
  systems with polytopic uncertainties.
\newblock {\em Int. J. Adapt. Control Signal Process.}, 20:291--304, 2006.

\bibitem{ZhHu10}
W.~Zhang, J.~Hu, and J.~Lian.
\newblock Quadratic optimal control of switched linear stochastic systems.
\newblock {\em Syst. Contr. Lett.}, 59:736--744, 2010.

\bibitem{Ti02}
P.H.E. Tiesinga.
\newblock Precision and reliability of periodically and quasiperiodically
  driven integrate-and-fire neurons.
\newblock {\em Phys. Rev. E}, 65(4):041913, 2002.

\bibitem{GrHo05}
T.~Griffin and S.~Hogan.
\newblock Dynamics of discontinuous systems with imperfections and noise.
\newblock In G.~Rega and F.~Vestroni, editors, {\em IUTAM Symposium on Chaotic
  Dynamics and Control of Systems and Processes in Mechanics.}, pages 275--285.
  Springer, 2005.

\bibitem{DiIo04}
M.F. Dimentberg and D.V. Iourtchenko.
\newblock Random vibrations with impacts: {A} review.
\newblock {\em Nonlinear Dyn.}, 36:229--254, 2004.

\bibitem{DiMe79}
M.F. Dimentberg and A.I. Menyailov.
\newblock Response of a single-mass vibroimpact system to white-noise random
  excitation.
\newblock {\em Z. Angew. Math. Mech.}, 59(12):709--716, 1979.

\bibitem{FoBr96}
M.~Fogli, P.~Bressolette, and P.~Bernard.
\newblock The dynamics of a stochastic oscillator with impacts.
\newblock {\em Eur. J. Mech. A-Solids}, 15(2):213--241, 1996.

\bibitem{SrPa05}
N~Sri~Namachchivaya and J.H. Park.
\newblock Stochastic dynamics of impact oscillators.
\newblock {\em J. Appl. Mech. Trans. ASME}, 72(6):862--870, 2005.

\bibitem{FeXu09b}
J.~Feng, W.~Xu, H.~Rong, and R.~Wang.
\newblock Stochastic responses of {D}uffing-{V}an der {P}ol vibro-impact system
  under additive and multiplicative random excitations.
\newblock {\em Int. J. Non-Linear Mech.}, 44:51--57, 2009.

\bibitem{ChEl05}
P.D. Christofides and N.H. El-Farra.
\newblock {\em Control of Nonlinear and Hybrid Process Systems. Designs for
  Uncertainty, Constraints and Time-Delays.}
\newblock Springer, New York, 2005.

\bibitem{FeZh06}
W.~Feng and J.-F. Zhang.
\newblock Stability analysis and stabilization control of multi-variable
  switched stochastic systems.
\newblock {\em Automatica}, 42:169--176, 2006.

\bibitem{MhEl05}
P.~Mhaskar, N.H. El-Farra, and P.D. Christofides.
\newblock Robust hybrid predictive control of nonlinear systems.
\newblock {\em Automatica}, 41:209--217, 2005.

\bibitem{LiSc00}
B.~Lindner and L.~Schimansky-Geier.
\newblock Coherence and stochastic resonance in a two-state system.
\newblock {\em Phys. Rev. E}, 61(6):6103--6110, 2000.

\bibitem{SiKu11}
D.J.W. Simpson and R.~Kuske.
\newblock Mixed-mode oscillations in a stochastic piecewise-linear system.
\newblock {\em Phys. D}, 240:1189--1198, 2011.

\bibitem{Fi88}
A.F. Filippov.
\newblock {\em Differential Equations with Discontinuous Righthand Sides.}
\newblock Kluwer Academic Publishers., Norwell, 1988.

\bibitem{WiDe00}
M.~Wiercigroch and B.~De~Kraker, editors.
\newblock {\em Applied Nonlinear Dynamics and Chaos of Mechanical Systems with
  Discontinuities.}, Singapore, 2000. World Scientific.

\bibitem{Br99}
B.~Brogliato.
\newblock {\em Nonsmooth Mechanics: Models, Dynamics and Control.}
\newblock Springer-Verlag, New York, 1999.

\bibitem{BlCz99}
B.~Blazejczyk-Okolewska, K.~Czolczynski, T.~Kapitaniak, and J.~Wojewoda.
\newblock {\em Chaotic Mechanics in Systems with Impacts and Friction}.
\newblock World Scientific, Singapore, 1999.

\bibitem{AwLa03}
J.~Awrejcewicz and C.~Lamarque.
\newblock {\em Bifurcation and Chaos in Nonsmooth Mechanical Systems.}
\newblock World Scientific, Singapore, 2003.

\bibitem{Ib09}
R.A. Ibrahim.
\newblock {\em Vibro-Impact Dynamics.}, volume~43 of {\em Lecture Notes in
  Applied and Computational Mechanics.}
\newblock Springer, New York, 2009.

\bibitem{Ts03}
C.K. Tse.
\newblock {\em Complex Behavior of Switching Power Converters.}
\newblock CRC Press, Boca Raton, FL, 2003.

\bibitem{Fi60}
A.F. Filippov.
\newblock Differential equations with discontinuous right-hand side.
\newblock {\em Mat. Sb.}, 51(93):99--128, 1960.
\newblock English transl. {\em Amer. Math. Soc. Transl.} 42(2):199--231, 1964.

\bibitem{CaGi06}
P.~Casini, O.~Giannini, and F.~Vestroni.
\newblock Experimental evidence of non-standard bifurcations in non-smooth
  oscillator dynamics.
\newblock {\em Nonlinear Dyn.}, 46(3):259--272, 2006.

\bibitem{LuGe06}
A.C.J. Luo and B.C. Gegg.
\newblock Stick and non-stick periodic motions in periodically forced
  oscillators with dry friction.
\newblock {\em J. Sound Vib.}, 291:132--168, 2006.

\bibitem{Jo03}
M.~Johansson.
\newblock {\em Piecewise Linear Control Systems.}, volume 284 of {\em Lecture
  Notes in Control and Information Sciences.}
\newblock Springer-Verlag, New York, 2003.

\bibitem{DiJo01}
M.~di~Bernardo, K.H. Johansson, and F.~Vasca.
\newblock Self-oscillations and sliding in relay feedback systems: Symmetry and
  bifurcations.
\newblock {\em Int J. Bifurcation Chaos}, 11(4):1121--1140, 2001.

\bibitem{Sc80}
Z.~Schuss.
\newblock {\em Theory and Applications of Stochastic Differential Equations.}
\newblock Wiley, New York, 1980.

\bibitem{Ga85}
C.W. Gardiner.
\newblock {\em Handbook of Stochastic Methods for Physics, Chemistry and the
  Natural Sciences.}
\newblock Springer-Verlag, New York, 1985.

\bibitem{GrVa99}
J.~Grasman and O.A. van Herwaarden.
\newblock {\em Asymptotic Methods for the Fokker-Planck Equation and the Exit
  Problem in Applications.}
\newblock Springer, New York, 1999.

\bibitem{Ts84}
Ya.Z. Tsypkin.
\newblock {\em Relay Control Systems.}
\newblock Cambridge University Press, New York, 1984.

\bibitem{FrPo02}
G.F. Franklin, J.D. Powell, and A.~Emami-Naeini.
\newblock {\em Feedback Control of Dynamic Systems.}
\newblock Prentice Hall, Upper Saddle River, NJ, 2002.

\bibitem{DoBi01}
R.C. Dorf and R.H. Bishop.
\newblock {\em Modern Control Systems.}
\newblock Prentice Hall, Upper Saddle River, NJ, 2001.

\bibitem{AsMu08}
K.J. {\AA}str\"{o}m and R.M. Murray.
\newblock {\em Feedback Systems. An Introduction for Scientists and Engineers.}
\newblock Princeton University Press, Princeton, NJ, 2008.

\bibitem{JoRa99}
K.H. Johansson, A.~Rantzer, and K.J. {\AA}str\"{o}m.
\newblock Fast switches in relay feedback systems.
\newblock {\em Automatica}, 35:539--552, 1999.

\bibitem{JoBa02}
K.H. Johansson, A.E. Barabanov, and K.J. {\AA}str\"{o}m.
\newblock Limit cycles with chattering in relay feedback systems.
\newblock {\em IEEE Trans. Automat. Contr.}, 47(9):1414--1423, 2002.

\bibitem{ZhFe10}
Y.~Zhao, J.~Feng, and C.K. Tse.
\newblock Discrete-time modeling and stability analysis of periodic orbits with
  sliding for switched linear systems.
\newblock {\em IEEE Trans. Circuits Systems I Fund. Theory Appl.},
  57(11):2948--2955, 2010.

\bibitem{TaOs09}
M.~Tanelli, G.~Osorio, M.~di~Bernardo, S.M. Savaresi, and A.~Astolfi.
\newblock Existence, stability and robustness analysis of limit cycles in
  hybrid anti-lock braking systems.
\newblock {\em Int. J. Contr.}, 82(4):659--678, 2009.

\bibitem{DiJo02}
M.~di~Bernardo, K.H. Johansson, U.~J\"{o}nsson, and F.~Vasca.
\newblock On the robustness of periodic solutions in relay feedback systems.
\newblock In {\em 15th Triennial World Congress, Barcelona, Spain}, 2002.

\bibitem{Ri84}
H.~Risken.
\newblock {\em The {F}okker-{P}lanck Equation: {M}ethods of Solution and
  Applications.}
\newblock Springer-Verlag, New York, 1984.

\bibitem{Sc10}
Z.~Schuss.
\newblock {\em Theory and Applications of Stochastic Processes.}
\newblock Springer, New York, 2010.

\bibitem{BeOr78}
C.M. Bender and S.A. Orszag.
\newblock {\em Advanced Mathematical Methods for Scientists and Engineers.}
\newblock International Series in Pure and Applied Mathematics. McGraw-Hill,
  New York, 1978.

\bibitem{KaSh84}
I.~Karatzas and S.E. Shreve.
\newblock Trivariate density of {B}rownian motion, its local and occupation
  times, with application to stochastic control.
\newblock {\em Ann. Prob.}, 12(3):819--828, 1984.

\bibitem{Gi60}
I.V. Girsanov.
\newblock On transforming a certain class of stochastic processes by absolutely
  continuous substitution of measures.
\newblock {\em Theory Prob. Appl.}, 5(3):285--301, 1960.

\bibitem{Ok03}
B.~{\O}ksendal.
\newblock {\em Stochastic Differential Equations: {A}n Introduction with
  Applications.}
\newblock Springer, New York, 2003.

\bibitem{KaSh91}
I.~Karatzas and S.E. Shreve.
\newblock {\em Brownian Motion and Stochastic Calculus.}
\newblock Springer, New York, 1991.

\bibitem{BeSh80}
V.E. Bene\u{s}, L.A. Shepp, and H.S. Witsenhausen.
\newblock Some solvable stochastic control problems.
\newblock {\em Stochastics}, 4:39--83, 1980.

\bibitem{QiZh02}
Z.~Qian and W.~Zheng.
\newblock Sharp bounds for transition probability densities of a class of
  diffusions.
\newblock {\em C.R. Acad. Sci. Paris, Ser. I}, 335:953--957, 2002.

\bibitem{QiRu03}
Z.~Qian, F.~Russo, and W.~Zheng.
\newblock Comparison theorem and estimates for transition probability densities
  of diffusion processes.
\newblock {\em Probab. Theory Relat. Fields.}, 127:388--406, 2003.

\bibitem{GrHe01}
M.~Gradinaru, S.~Herrmann, and B.~Roynette.
\newblock A singular large deviations phenomenon.
\newblock {\em Ann. I. H. Poincar\'{e}}, 37(5):555--580, 2001.

\bibitem{Zh90}
W.~Zhang.
\newblock Transition density of one-dimensional diffusion with discontinuous
  drift.
\newblock {\em IEEE Trans. Automat. Contr.}, 35(8):980--985, 1990.

\bibitem{CoDi12}
A.~Colombo, M.~di~Bernardo, S.J. Hogan, and M.R. Jeffrey.
\newblock Bifurcations of piecewise smooth flows: perspectives, methodologies
  and open problems.
\newblock Submitted to: {\em Phys. D}, 2012.

\end{thebibliography}

\end{document}